\documentclass[a4paper, 12pt]{article}

\RequirePackage{amsthm,amsmath,amsfonts,amssymb}
\RequirePackage[round]{natbib}
\RequirePackage[colorlinks,citecolor=black,urlcolor=black]{hyperref}
\RequirePackage{graphicx}

\usepackage{paralist}
\usepackage{setspace}

\usepackage[linesnumbered, ruled,vlined]{algorithm2e}
\usepackage{textcomp}
\usepackage{mathtools}
\usepackage{enumitem} 
\usepackage{mleftright}
\usepackage{multirow}
\usepackage[figuresright]{rotating}
\usepackage[utf8]{inputenc}
\usepackage[english]{babel}
\usepackage{filecontents}
\usepackage{sectsty}
\usepackage{titletoc}
\usepackage[toc,page]{appendix}
\graphicspath{ {figures/} }
\usepackage{array}
\usepackage{float}
\usepackage{scalerel,stackengine}
\usepackage{accents}
\usepackage{color, graphics, graphicx}
\usepackage{xcolor,colortbl}
\usepackage{undertilde}
\usepackage[T1]{fontenc}
\usepackage{lmodern}
\usepackage{hyperref}
\newcommand\tenq[2][1]{%
\def\useanchorwidth{T}%
\ifnum#1>1%
\stackunder[0pt]{\tenq[\numexpr#1-1\relax]{#2}}{\scriptscriptstyle\thicksim}%
\else%
\stackunder[1pt]{#2}{\scriptscriptstyle\thicksim}%
\fi%
}

\usepackage{subfigure}
\numberwithin{equation}{section}

\newtheorem{theorem}{Theorem}[section]
\newtheorem{lemma}[theorem]{Lemma}
\newtheorem{proposition}[theorem]{Proposition}

\theoremstyle{remark}
\newtheorem{assumption}{Assumption}

\numberwithin{theorem}{section}

\newcommand{\pr}{^{\prime}}
\newcommand{\n}{^{(n)}}
\newcommand{\bth}{{\boldsymbol \theta}}
\newcommand{\tbth}{{\scriptstyle{\boldsymbol{\theta}}}}

\DeclareMathAlphabet\mathbfcal{OMS}{cmsy}{b}{n}

\newcommand{\Zb}{{\bf Z} }
\newcommand{\mub}{{\pmb\mu} }

\usepackage{stackengine}
\stackMath

\def\pms{\mspace{-1mu}{\scriptscriptstyle\pm}}

\def\s{\sum_{t=1}^n}

\def\bth{\mbox{\boldmath$\theta$}}

\def\bLam{\mbox{\boldmath$\Lambda$}}

\def\btau{\mbox{\boldmath$\tau$}}
\def\bDelta{\mbox{\boldmath$\Delta$}}

\def\bSigma{\mbox{\boldmath$\Sigma$}}
\def\bXi{\mbox{\boldmath$\Xi$}}
\def\bmu{\mbox{\boldmath$\mu$}}

\def\bGamma{\mbox{\boldmath$\Gamma$}}
\def\bOmega{\mbox{\boldmath$\Omega$}}
\def\bUpsilon{\mbox{\boldmath$\Upsilon$}}
\def\bvp{\mbox{\boldmath$\varphi$}}
\def\bepsilon{\mbox{\boldmath$\epsilon$}}

\def\vp{\varphi}

\def\vp{\varphi}

\def\R{{\mathbb R}}
\def\Z{{\mathbb Z}}
\def\D{{\bf D}}
\def\N{{\mathbb N}}

\def\K{{\bf K}}
\def\E{{\bf E}}
\def\J{{\bf J}}

\def\u{{\bf u}}

\def\M{{\bf M}}
\def\bz{{\bf z}}

\def\0{{\bf 0}}
\def\W{{\bf W}}

\def\Q{{\bf Q}}
\def\S{{\bf S}}

\def\A{{\bf A}}
\def\B{{\bf B}}

\def\C{{\bf C}}
\def\DD{{\bf D}}
\def\F{{\bf F}}
\def\H{{\bf H}}
\def\hbf{{\bf h}}
\def\G{{\bf G}}
\def\I{{\bf I}}
\def\X{{\bf X}}
\def\x{{\bf x}}
\def\y{{\bf y}}

\def\ZZ{{\bf Z}}

\def\0{{\bf 0}}

\def\c{{\bf c}}

\def\e{{\bf e}}

\def\1{{\bf 1}}

\newcommand{\tetb}{{\pmb \theta} }

\title{Center-outward Rank- and Sign-based VARMA Portmanteau Tests:\\ Chitturi, Hosking, and Li--McLeod revisited
} 
\date{} 
\author{
Marc Hallin\footnote{Corresponding author: Marc Hallin, Universit\'e Libre de Bruxelles CP 114/4, Avenue Franklin Roosevelt 50,  
B-1050 Bruxelles, Belgium. }\\
{\small ECARES and D\'epartement de Math\'ematique}\\
{\small  Universit\'e Libre de Bruxelles,  Brussels,  
Belgium} \\ 
{\small E-mail: mhallin@ulb.ac.be}\vspace{3mm}
\\  \, {\small and}  \vspace{3mm}\\\ 
Hang Liu\\
{\small International Institute of Finance, School of Management}\\ {\small University of Science and Technology of China}\\
{\small Hefei, 
 Anhui,   
 China} \\ 
{\small E-mail: hliu01@ustc.edu.cn}
}

\begin{document}
\maketitle

\begin{abstract}
{\small{The pseudo-Gaussian portmanteau  tests of  Chitturi, Hosking, and Li and McLeod for VARMA models are revisited from a Le Cam perspective,  providing a  precise and  more rigorous description of the asymptotic behavior of the multivariate portmanteau  test statistic, which depends on the dimension $d$ of the observations, the number $m$ of lags involved, and the length $n$ of the observation period. Then, based on the concepts of center-outward ranks and signs recently developed (Hallin, del Barrio, Cuesta-Albertos, and Matr\' an, {\it Annals of  Statistics} 49,  1139--1165, 2021),  a class of multivariate  rank- and sign-based portmanteau test statistics is proposed which, under the null hypothesis and under a broad family of innovation densities,  can be approximated by an asymptotically chi-square variable. The asymptotic properties of these tests  are derived;  simulations demonstrate their advantages  over their classical pseudo-Gaussian counterpart. } }
\end{abstract}


\section{Introduction}

\subsection{The Gaussian multivariate portmanteau test}

The so-called {\it portmanteau test}  certainly ranks among the most popular and most widely used testing procedures in time series analysis. It is simple, intuitive, apparently well understood, and naturally complements eye-inspection of residual correlograms. 

In their univariate forms, portmanteau test statistics were first introduced by  \cite{BP70} as a sum of  squared residual autocorrelations of orders one through $m$ associated with the estimation (usually, a least square one) of the model parameters. A later version by \cite{jBox78}  is taking into account the fact that a residual autocorrelation of order $k$, when computed from a series of length $n$, is based on a sum of $(n-k)$ terms only;  while improving finite-sample performance, that modification has no asymptotic impact, though.

\cite{Chit74} for the VAR case, \cite{Hosk1980}, and \cite{LiMcL81} for the VARMA case extended the Ljung-Box-Pierce test to the multivariate context by replacing sums of squared residual autocorrelations with  sums of normalized squared elements of (estimated)  residual cross-covariance matrices; modified versions in the spirit of \cite{jBox78}  also are proposed by these authors.

A widespread opinion is that the null distribution of the multivariate portmanteau statistic, in the $d$-dimensional VARMA$(p,q)$ case, is asymptotically chi-square with $d^2(m-p-q)$ degrees of freedom under a ``broad'' class of innovation densities, where $m$ is the number of lags in the test statistic.  This is, actually, the statement in Theorem~2 of \cite{Hosk1980}. In this statement, Hosking is not very precise about what is to be understood with ``asymptotically chi-square.'' Actually, if only the series length $n$ goes to infinity, the claim is incorrect.  In his  Appendix, he recommends   $m= O(n^{1/2})$; but this clearly precludes an asymptotic distribution with finitely many degrees of freedom. \cite{Chit74}, and \cite{LiMcL81} are more cautious, with a somewhat vague  claim that this asymptotic null distribution is ``approximately chi-square with $d^2(m-p-q)$ degrees of freedom for $m$ and $n$ (the series length) large enough.'' 

As we shall see, these statements, at best, are imprecise and  the asymptotic distribution of the portmanteau test statistic under the null is not chi-square with $d^2(m-p-q)$ degrees of freedom---not even under Gaussian innovation densities.  The chi-square  critical values resulting from these statements, nevertheless,  are routinely used in the daily  implementation of the  test.

In the univariate case, a similar problem was detected by 
Taniguchi and Amano~(2010), who show that the classical univariate portmanteau test statistics of \cite{BP70} and \cite{jBox78} 
 are not asymptotically chi-square if the number $m$ of lagged residuals in the test statistic is finite. 
 Based on this observation, they also propose a modified Whittle likelihood ratio test which is asymptotically chi-square. In a more general setting of single-output regression with ARMA errors, \cite{Akashi2018, Akashi2021} propose a likelihood ratio-based portmanteau test incorporating the Whittle likelihood ratio test of \cite{Taniguchi2010} as a special case and provide a sufficient condition, in terms of   Fisher information, under which their test statistic is asymptotically chi-square.

Half a century after its introduction, thus, one of the most  popular tests in time series analysis still relies on   vague or even faulty asymptotic statements. The first objective of this paper is to fix these asymptotic results. Rather than modifying  the classical test statistics (as \cite{Taniguchi2010} and \cite{Akashi2018, Akashi2021} are doing in the univariate case), we are providing a   precise account and a rigorous derivation of the asymptotic behavior of the portmanteau test statistic  for VARMA models. Taking advantage of this, we then propose, based on the measure-transportation-based concepts of center-outward ranks and signs recently developed in \cite{Hallinetal2021} (see \cite{H22} for a nontechnical survey) a class of   rank-based portmanteau tests and establish their asymptotic distributions under the null. 
Simulations demonstrate the advantages of these tests, in terms of power and resistance to outliers, over  the classical Chitturi-Hosking-Li-McLeod pseudo-Gaussian ones under non-Gaussian innovation densities. 

Throughout, the theoretical tools we are using are borrowed from Le Cam's powerful asymptotic theory of statistical experiment---a theory that was not available,  fifty years ago,  to Chitturi, Hosking, Li, and McLeod, who instead are making intensive use of Taylor expansions.

\subsection{Center-outward rank-based multivariate portmanteau test}

Despite their popularity, the Gaussian portmanteau tests have some undesirable drawbacks. First, their validity and consistency rely on the assumption of  finite fourth-order innovation moments, an assumption that may not hold in economic and financial practice, where data are usually heavy-tailed. Second, the finite-sample performance of pseudo-Gaussian tests often is quite poor under non-Gaussian innovations (see  Section~1.2 in~\cite{HLL2022} for numerical examples). Last but not least, the pointwise chi-square approximation of the null distribution of the pseudo-Gaussian portmanteau statistic under given  innovation density~$f$ is far from uniform with respect to $f$. As a consequence, the $\sup_f$ of the size under innovation density $f$ of pseudo-Gaussian tests, in general, does not converge to the nominal asymptotic level~$\alpha$  (see, e.g.,  Section~1.1 in~\cite{HHH22}): as a semiparametric test, thus, the pseudo-Gaussian portmanteau test fails to satisfy the asymptotic probability level condition.

Thanks to distribution-freeness, classical (univariate) rank-based tests  are escaping that asymptotic size problem: for the ARMA case, see, e.g., \cite{HP88,HP94}.   In a multivariate or multiple-output context, however,   due to the fact that   no canonical ordering is available  in~$\mathbb{R}^d $ for~$d\geq~\!2$, a long-standing problem has been: ``what are ranks and signs in dimension $d\geq 2$?'' Various notions of ranks and signs have been proposed in the literature, including the componentwise ranks \citep{Puri1971}, the spatial ranks \citep{Oja2010}, the depth-based ranks \citep{Liu1992}, and the  Mahalanobis ranks and signs \citep{HP04}. None of these ranks or signs are distribution-free under the whole family of absolutely continuous distributions, though. The {\it center-outward ranks and signs} recently proposed by \cite{Chernozhukov2017} and \cite{Hallinetal2021} are not only  distribution-free  under  absolutely continuous distributions; they also are maximal ancillary (see \cite{Hallinetal2021}). They  have been applied quite successfully to various statistical models of daily statistical importance, including  multiple-output regression and MANOVA  \citep{HHH22}, goodness-of-fit  \citep{GhS22}, test of  vector independence  \citep{DebSen2021, Shietal22}, multivariate quantile regression\linebreak \citep{delBarrio2022}, R-estimation for VARMA models \citep{HLL2022} and rank-based order selection of VAR models \citep{HLL2023}. 

In this paper, we propose a class of asymptotically distribution-free
 portmanteau tests for VARMA models based on   residual center-outward ranks and signs. These residuals are the estimated residuals based on an estimation~$\tenq{\tetb}\n$ of the parameter $\tetb$ of the null VARMA model: the resulting ranks and signs, thus, are {\it aligned} ones, failing to achieve exact finite-sample distribution-freeness.  For  sufficiently large $m$, however, the test statistic is  arbitrarily close to an  asymptotically chi-square ($d^2(m-p-q)$ degrees of freedom) oracle statistic based on the ranks and signs of the {\it exact residuals}---that is, based on the actual $\tetb$ value, as opposed to the {\it aligned} ones; the ranks and signs of these exact residuals are fully distribution-free but, of course, cannot be computed from the observations.

 While the pseudo-Gaussian portmanteau test  requires plugging-in  the Gaussian quasi-likelihood VARMA estimator $\hat\tetb\n$ for the unspecified VARMA parameter $\tetb$, our rank-based tests rely on plugging-in  the center-outward R-estimators $\tenq{\tetb}\n$---as derived in \cite{HLL2022}---based on the same score function as the test statistic itself. 
  A numerical study reveals that, when based on Gaussian scores, our rank-based test outperforms the classical pseudo-Gaussian one  under a broad range of innovation distributions. This is probably the sign of a general Chernoff-Savage property \citep{ChernoffS} for which, however, we have no proof in this context.

The paper is organized as follows. Section~\ref{Sec:notation} introduces the VARMA model and the local asymptotic normality of the model under some mild regularity assumptions. Section~\ref{Sec:GaussTest} revisits the Gaussian multivariate portmanteau tests of \cite{Chit74}, \cite{Hosk1980} and \cite{LiMcL81} and re-establishes their asymptotic properties  under precise form  via  Le Cam's asymptotic theory. Section~\ref{Sec.RankPortmanteau} proposes a class of  center-outward rank-based portmanteau test statistics and similarly establishes their asymptotic properties. A brief numerical analysis of the finite-sample size and power of these tests is conducted in Section~\ref{Sec:Numerical}. Section~6 concludes. All proofs   are collected in the~Appendix.


%


\section{Notation and general setting}\label{Sec:notation}

In this section, we introduce the VARMA($p, q$) model (Section~\ref{Sec:VARMA}) and the local asymptotic normality property of the model (Section~\ref{secLAN}), which is the essential tool in our analysis of the asymptotic behavior of the Gaussian and center-outward rank-based multivariate portmanteau test statistics.

\subsection{VARMA models}\label{Sec:VARMA}

Consider the $d$-dimensional VARMA($p, q$) model
\begin{equation}
\Big(\I_d - \sum_{i = 1}^p \A_i L^i \Big) \X_t = \Big(\I_d + \sum_{j = 1}^q \B_j L^j\Big) \bepsilon_t, \quad t \in \Z, \vspace{-2mm}\label{VARMA_mod}
\end{equation}
where $\A_1, \ldots , \A_p, \B_1, \ldots , \B_q$ are $d \times d$ matrices, $L$ denotes the lag operator, and~$\{\bepsilon_t; t \in \Z\}$ is some i.i.d.\ mean-zero $d$-dimensional  white noise process  with density $f$. Denoting by 
$$\bth := \big((\text{vec}({\A_1}))^\prime, \ldots , (\text{vec}({\A_p}))^\prime, (\text{vec}({\B_1}))^\prime, \ldots , (\text{vec}({\B_q}))^\prime\big)^\prime,$$  
 (where $^\prime$ indicates transposition) the $(p+q)d^2$-dimensional VARMA parameter, we throughout   assume that $\bth$ satisfies the following very classical conditions ensuring  identifiability and the existence of a stationary and invertible solution to~\eqref{VARMA_mod}.
 
 \begin{assumption}\label{ass.stationary}
\begin{enumerate}
\item[{\it (i)}] All solutions of the determinantal equations
$$\text{det$\left(\I_d - \sum_{i = 1}^p \A_i z^i \right) = 0\quad$ and \quad det$\left( \I_d + \sum_{i = 1}^q \B_i z^i \right) = 0$, \quad$z \in \mathbb{C}$}$$
 lie outside the unit ball in $\mathbb{C}$;
\item[{\it (ii)}]  det$(\A_p) \neq 0 \neq$det$(\B_q)$; 
\item[{\it (iii)}]   
$\I_d - \sum_{i = 1}^p \A_i z^i$ and $\I_d + \sum_{i = 1}^q \B_i z^i$ have no common left
factors other than
$\I_d$.  
\end{enumerate}
\end{assumption}
\noindent Denote by $\boldsymbol\Theta_{p,q}$ the VARMA($p, q$)   {\it parameter space}---namely, the set of   all $\bth$ values satisfying Assumption~\ref{ass.stationary}: $\boldsymbol\Theta_{p,q}$ is a finite collection of open connected subsets of ${\R}^{(p+q)d^2}$. 


%

Let  $\X^{(n)}\coloneqq \{\X^{(n)}_1, \ldots , \X^{(n)}_n\}$ (superscript\! $^{(n)}\!$  omitted whenever possible) be an observed  finite realization of some solution of \eqref{VARMA_mod}. For any~$\bth\in~\!\boldsymbol\Theta_{p,q}$,  this  observation 
 $\X^{(n)}$
   is asymptotically (as $n\to\infty$) stationary and invertible; associated with  
    any~$(p+q)$-tuple~$(\X^{(n)}_0, \ldots , \X^{(n)}_{-p+1}, \bepsilon^{(n)}_1, \ldots , \bepsilon^{(n)}_{-q+1})$  of initial values, it determines, for any $\bth$, an $n$-tuple 
    of {\it residuals} 
  \[\Zb^{(n)}_t(\bth)\coloneqq \Big(\I_d - \sum_{i = 1}^p \A_i L^i \Big) \X\n_t  -  \sum_{j = 1}^q \B_j L^j \bepsilon_t, \quad t=1,\ldots,n,
  \]
 which  can be computed recursively.  
These residuals are i.i.d.\ with density~$f$ and  $\Zb^{(n)}_t(\bth)$ coincides with $\bepsilon_t$ for any $t$ iff $\X^{(n)}$ is a solution of \eqref{VARMA_mod} with  parameter value $\bth$.

%

\subsection{Local asymptotic normality (LAN)}\label{secLAN}

Our asymptotic analysis relies on the local asymptotic normality (LAN) property of VARMA models. That property requires mild regularity conditions on the density $f$ of $\bepsilon_t$ . 
 
\begin{assumption}\label{ass.den}{\rm 
\begin{enumerate}
\item[{\it (i)}] The innovation density $f$ belongs to the class   of {\it non-vanishing} Lebesgue densities on~$\R^d$, i.e., $f(\x)>0$ for all $\x\in{\mathbb R}^d$;
\item[{\it (ii)}]  $\displaystyle{\int \x f(\x) \mathrm{d}\mu = \0}$ and $\displaystyle{ \int \x \x^\prime f(\x) \mathrm{d}\mu = \bXi}$ where~$\bXi$ is finite and positive definite; 
\item[{\it (iii)}] $f^{1/2}$ is {\it mean-square differentiable} with {\it mean-square gradient}~$\DD f^{1/2}$, that is, 
there exists a square-integrable   vector $\DD f^{1/2}$ such that, for any sequence~$\hbf\in\mathbb{R}^d$ such  that~$\0 \neq \hbf \rightarrow \0$,
$${(\hbf^\prime \hbf)^{-1} \int \left[ f^{1/2}(\x + \hbf) - f^{1/2}(\x) - \hbf^\prime \DD  f^{1/2} (\x)\right]^2 \mathrm{d}\mu \rightarrow 0;}$$
\item[{\it (iv)}] letting 
$\bvp_{f}(\x) \coloneqq (\vp_{f1}(\x), \ldots , \vp_{fd}(\x))^\prime \coloneqq -2 \DD f^{1/2}(\x)/f^{1/2}(\x)$ (the location score function), 
 $\int [\vp_i (\x)]^4 f(\x) \mathrm{d}\mu < \infty,\ i = 1, \ldots , d$;
\item[{\it (v)}]   the   function $\x\mapsto\bvp_f(\x)$ is piecewise Lipschitz, i.e., there exists $K\in\R$ and a finite measurable partition of $\R^d$ into $J$ non-overlapping subsets~$I_1,  \ldots , I_J$ such that,  for all $\x, \y$ in~$I_j$, $j = 1, \ldots , J$,
\[
\Vert \bvp_f(\x) - \bvp_f(\y)\Vert  \leq K \Vert \x - \y\Vert .\]
\end{enumerate}}
\end{assumption}
\noindent Conditions {\it (i)} and {\it (iii)} are standard in the context of LAN; {\it (ii)} and {\it (iv)}  guarantee the existence of a full-rank information matrix for $\bth$. Under {\it (v)}, the impact of initial values is  asymptotically negligible in mean square norm; without loss of generality, thus, we henceforth assume that $$\X^{(n)}_0= \ldots = \X^{(n)}_{-p+1}=\0= \bepsilon^{(n)}_1=\ldots = \bepsilon^{(n)}_{-q+1}.$$  
Denote by $\mathcal{F}_d$ the class of densities $f$ satisfying Assumption~\ref{ass.den}.

Letting ${\rm P}^{(n)}_{\tbth ;f}$ denote the distribution of~$\X^{(n)}$ under parameter value~$\bth$ and innovation density~$f$, write
$$L^{(n)}_{\tbth + n^{-1/2}\btau^{(n)}/\tbth; f}\coloneqq \log\frac{ {\rm d}{\rm P}^{(n)}_{\tbth + n^{-1/2}\btau^{(n)};f} }{ {\rm d}{\rm P}^{(n)}_{\tbth;f}}
,$$   where~$\btau^{(n)}$ is a bounded sequence of $\mathbb{R}^{(p+q) d^2}$, for the log-likelihood ratio of~${\rm P}^{(n)}_{\tbth + n^{-1/2}\btau^{(n)};f}$ with respect to ${\rm P}^{(n)}_{\tbth ;f}$ computed at $\X\n$.  \smallskip

Define 
 \begin{align}\label{Sf1}
\bGamma_{f}^{(n)}(\bth) \coloneqq   n^{-1/2}  \left((n-1)^{1/2} \big(\text{vec}({\bGamma_{1, f}^{(n)}(\bth)})\big)^\prime, \ldots,\right. &\nonumber \\
 (n-i)^{1/2} \big(\text{vec}({\bGamma_{i; f}^{(n)}(\bth))}\big)^\prime, \ldots&\left.  , \big(\text{vec}({\bGamma_{n-1, f}^{(n)}(\bth)})\big)^\prime\right)^\prime,
\end{align}
with the so-called  $f$-{\it cross-covariance matrices}
\begin{equation}\label{Gamma1}
\bGamma_{i, f}^{(n)}(\bth) \coloneqq (n-i)^{-1} \sum_{t=i+1}^n \bvp_{f}(\ZZ_t^{(n)}(\bth)) \ZZ_{t-i}^{(n)\prime} (\bth)\quad i=1,\ldots, n-1.
\end{equation}
 Let $\C\n_{\tbth} \coloneqq \left(\c_{1,\tbth}, \ldots, \c_{n-1,\tbth}\right)$ with
 \begin{equation}\label{def.ci}
\c_{i,\tbth} 
\coloneqq 
\begin{pmatrix}
  \sum_{j=0}^{i-1} \sum_{k=0}^{\min (q, i-j-1)} (\G_{i-j-k-1} \B_k) \otimes \H_j^\prime  \\
\vdots \\
  \sum_{j=0}^{i-p} \sum_{k=0}^{\min (q, i-j-p)} (\G_{i-j-k-p} \B_k) \otimes \H_j^\prime  \\
  \I_d \otimes \H_{i-1}^\prime  \\
\vdots \\
  \I_d \otimes \H_{i-q}^\prime 
\end{pmatrix} 
,\  i = 1, \ldots, n-1,
\end{equation}
where $\G_u$ and $\H_u$, $u\in \mathbb{Z}$ are the {\it Green matrices} associated with the autoregressive and moving average 
 operators, respectively, in \eqref{VARMA_mod}---namely, the matrix coefficients of the inverted linear difference operators~$\big(\A (L)\big)^{-1}$ and~$\big(\B (L)\big)^{-1}\!$:
$$\sum_{u=0}^\infty \G_u z^u\! = \left(  \I_d - \sum_{i=1}^p  \A_i z^i\right)^{-1}\text{and}\quad\!\sum_{u=0}^\infty \H_u z^u\! = \left( \I_d + \sum_{i=1}^q \B_i z^i\right)^{-1}\!\!\! ,  \, z \in \mathbb{C}, |z| < 1.$$

More constructive (recursive) definitions of Green's matrices can be found in Section~1.1 of \cite{H86} and Section~3 of \cite{GH95}  but are not needed here, and we skip them for the sake of space: the only thing  we need  to recall is the fact that,  under Assumption~\ref{ass.stationary},  $\Vert \G_u\Vert $ and~$\Vert \H_u\Vert$, just as the moduli of all solutions of the homogeneous difference equations associated with~$\A (L)$ and~$\B (L)$,   are  exponentially  decreasing  as $u\to\infty$, which implies that  $\Vert\c_{i,\tbth}\Vert$  exponentially   decreases as $i\to\infty$. Defining 
\begin{align}
\bDelta^{(n)}_{f} (\bth) &\coloneqq   \sum_{i=1}^{n-1} \c_{i,\tbth}  (n-i)^{1/2} {\rm vec}({\bGamma_{i, f}^{(n)}(\bth)}) \label{Delta}   = n^{1/2} \C\n_{\tbth} \bGamma_{f}^{(n)}(\bth), 
\end{align}
(the central sequences), the following LAN result is established in~\cite{GH95}, Proposition~3.1.

\begin{proposition}[Garel and Hallin (1995)]\label{Prop.LAN1}
Let Assumptions~\ref{ass.stationary}  and \ref{ass.den} hold. Then, for  any bounded sequence~$\btau^{(n)}$ in~$\R^{(p+q) d^2}\!$, under ${\rm P}^{(n)}_{\tbth;f}$, as $n \rightarrow \infty$, 
\begin{equation}\label{lik1}
L^{(n)}_{\tbth + n^{-1/2}\btau^{(n)}/\tbth; f} = \btau^{{(n)}\prime} \bDelta^{(n)}_{f}(\bth) - \frac{1}{2} \btau^{(n)\prime} \bLam_{f} (\bth) \btau^{(n)} + o_{\rm P}(1),
\end{equation}
with  $\bDelta^{(n)}_{f} (\bth) $  (the central sequence)   defined in \eqref{Delta} and $(p+q) d^2\times (p+q) d^2$ symmetric and positive definite $\bLam_{f} (\bth)$ (the information matrix---see equation~(3.16) in \cite{GH95} for an explicit form). Still under~${\rm P}^{(n)}_{\tbth;f}$,~$\bDelta^{(n)}_{f} (\bth) $ is asymptotically normal with mean $\0$ and   covariance  $\bLam_{f} (\bth)$. 
\end{proposition}

\section{The Chitturi-Hosking-Li-McLeod multivariate portmanteau test}\label{Sec:GaussTest}

In this section, we are revisiting the pseudo-Gaussian multivariate portmanteau test of Chitturi-Hosking-Li-McLeod and,  adopting a  Le Cam approach, clarify their asymptotic results. 

\subsection{Residual cross-covariance matrices}\label{sec.GaussCovMat}

For a Gaussian density $f$, the score function is $\bvp_f(\bz) = -{\bf \Sigma}^{-1}\bz$, where ${\bf \Sigma}$ is the covariance matrix of $\bepsilon_t$, and the  $f$-cross-covariance matrices \eqref{Gamma1}  reduce~to \vspace{-2mm}
\begin{equation}\label{gammadef}
{\bGamma}_{i; \mathcal{N}}^{(n)}(\bth) \coloneqq - (n-i)^{-1} {\bf \Sigma}^{-1} \sum_{t=i+1}^n   \ZZ_{t}(\bth) \ZZ^{\prime}_{t-i}(\bth), \quad i = 1, \ldots , n - 1.
\vspace{-2mm}\end{equation}
Note that ${\bGamma}_{i; \mathcal{N}}^{(n)}(\bth) $ differs from the traditional lag-$i$ cross-covariance matrix~$(n-i)^{-1}  \sum_{t=i+1}^n   \ZZ_{t}(\bth) \ZZ^{\prime}_{t-i}(\bth)$  by a left factor $-{\bf \Sigma}$. 

Proposition~\ref{asy.Gami.N} provides the asymptotic distribution of~${\bGamma}_{i; \mathcal{N}}^{(n)}(\bth)$ under~${\rm P}^{(n)}_{\tbth ;f}$ and contiguous sequences of alternatives ${\rm P}^{(n)}_{\tbth + n^{-1/2}\btau ;f}$. See  the Appendix for a proof. 
 In practice, of course, $ {\bf \Sigma}$ remains unspecified and has to be estimated by some $ {\bf \Sigma}\n$---typically, the residual empirical covariance matrix. As long as~$ {\bf \Sigma}\n$  is consistent, substituting it for $ {\bf \Sigma}$ in \eqref{gammadef}, in view of Slutsky's Lemma, has no impact on  asymptotic distributions. Rather than introducing
cumbersome additional notation, thus, we pursue with the $ {\bf \Sigma}$-based definition of ${\bGamma}_{i; \mathcal{N}}^{(n)}(\bth)$. As for its rank-based counterparts defined in Section~\ref{SecSR}, 
 they do not involve any $\boldsymbol\Sigma$ nor any of the parameters of the actual innovation distribution.

\begin{proposition}\label{asy.Gami.N}
Let Assumptions~\ref{ass.stationary} and \ref{ass.den} hold. Then, 
 for any positive integers $i\neq j$, the vectors
 $$(n-i)^{1/2} {\rm vec}({\bGamma}_{i; \mathcal{N}}^{(n)}(\bth)) \quad\text{and}\quad (n-j)^{1/2} {\rm vec} ({\bGamma}_{j; \mathcal{N}}^{(n)}(\bth))$$
  are jointly asymptotically normal, 
with mean $(\0\pr, \0\pr)\pr$ under  ${\rm P}^{(n)}_{\tbth ;f}$,  mean\vspace{-1mm}
\begin{equation}\label{aslinN}
\left(({\bf \Sigma}\otimes {\bf \Sigma}^{-1}) \c_{i,\tbth}\pr  \btau)\pr, (({\bf \Sigma}\otimes {\bf \Sigma}^{-1})\c_{j,\tbth}\pr \btau)\pr\right)\pr
\end{equation}
under ${\rm P}^{(n)}_{\tbth + n^{-1/2}\btau; f}$, and   covariance  $\left(
\begin{array}{cc} {\bf \Sigma}\otimes {\bf \Sigma}^{-1} & \0 \\ \0& {\bf \Sigma}\otimes {\bf \Sigma}^{-1} \end{array}
\right)$  under both.
\end{proposition}

To determine the asymptotic behavior of the   Gaussian portmanteau test, the asymptotic linearity of $(n-i)^{1/2} {\rm vec}({\bGamma}_{i; \mathcal{N}}^{(n)}(\bth))$ is required. That result is established in Lemma~4 of \cite{HP05} under the quite restrictive assumption that  the distribution of $\bepsilon_t$ is elliptically symmetric. Here,  we make the assumption that asymptotic linearity  holds for  general $f$ satisfying Assumption~\ref{ass.den};  the form of the linear term on the right-hand side of~\eqref{eq.linN} below follows from the form of the shift matrices in~\eqref{aslinN}---themselves a consequence of Le Cam's third Lemma.

\begin{assumption}\label{asylinN}
For any positive integer $i$,  \vspace{-2mm}
\begin{align}\label{eq.linN}
(n-i)^{1/2}\left[{\rm vec} ({\bGamma}^{(n)}_{i; \mathcal{N}}(\bth + n^{-1/2}\btau)) \right.&\left. - {\rm vec} ({\bGamma}^{(n)}_{i; \mathcal{N}}(\bth)) \right] = ({\bf \Sigma}\otimes {\bf \Sigma}^{-1}) \c_{i,\tbth}\pr \btau + o_{\rm P}(1) 
\end{align}\vspace{-9mm}

\noindent under  ${\rm P}^{(n)}_{\tbth ;f}$, as $n\to\infty$.
\end{assumption} 
\noindent We refer to Section~4 of \cite{vdAetal2015} for primitive conditions.

The Gaussian portmanteau test is based on the plug-in of $\hat{\bth}\n_{\cal N}$, the Gaussian quasi-maximum likelihood estimator (QMLE) of $\bth$, in the test statistic to be defined. Recall that $\hat{\bth}\n_{\cal N}$ is defined as the solution of 
 $\bDelta\n_{\cal N}(\bth) = \0$ 
where 
$$\bDelta\n_{\cal N}(\bth) \coloneqq  \sum_{i=1}^{n-1} \c_{i,\tbth}    (n-i)^{1/2} {\rm vec}({\bGamma_{i; \mathcal{N}}^{(n)}(\bth)})$$
is the Gaussian central sequence. A finite fourth-order moment assumption on~$\bepsilon_t$ is required in order for  $\hat{\bth}\n_{\cal N}$ to be root-${n}$-consistent; see, e.g., M\' elard~(2022) for a recent proof.\medskip

\begin{assumption}\label{ass.fourMom}
${\rm E} \left[  {\rm vec}(\bepsilon_t \bepsilon\pr_t)\left( {\rm vec}(\bepsilon_t \bepsilon\pr_t)\right)\pr\right] < \infty.$
\medskip
\end{assumption}

In view of \eqref{eq.linN}, the plug-in impact of substituting $\hat{\bth}\n_{\cal N}$ for the actual $\bth$ in a quadratic form of~$(n-i)^{1/2}{\rm vec} ({\bGamma}^{(n)}_{i; \mathcal{N}}(\bth))$, $i=1,\ldots, m$  is asymptotically non-negligible. 
 Taking into account the fact that $\bDelta\n_{\cal N}(\hat{\bth}\n_{\cal N}) = \0$, however, that impact can be neutralized in view of Lemma~\ref{Lemport} below.
\color{black}

Denote by  ${\bGamma}_{i; {\cal N}}^{(n)*}(\bth)$ the matrix of residuals in the regression of  ${\bGamma}_{i; {\cal N}}^{(n)}(\bth)$'s entries with respect to ${\bDelta}_{{\cal N}}^{(n)}(\bth)$ in their joint asymptotic (under  ${\rm P}^{(n)}_{\tbth ;f}$) distribution; namely, let     
\begin{align}
&(n-i)^{1/2} {\rm vec} ({\bGamma}^{(n)*}_{i; {\cal N}}(\bth)) \nonumber \\
&\ \coloneqq (n-i)^{1/2} {\rm vec} ({\bGamma}^{(n)}_{i; {\cal N}}(\bth)) 
- ({\bf \Sigma}\otimes {\bf \Sigma}^{-1}) \c_{i,\tbth}\pr\Big( \sum_{i=1}^{n-1} \c_{i,\tbth} ({\bf \Sigma}\otimes {\bf \Sigma}^{-1}) \c_{i,\tbth}\pr\Big)^{-1} {\bDelta}^{(n)}_{{\cal N}}(\bth) \label{=*} \\
&\ = (n-i)^{1/2} {\rm vec} ({\bGamma}^{(n)}_{i; {\cal N}}(\bth)) - ({\bf \Sigma}\otimes {\bf \Sigma}^{-1}) \c_{i,\tbth}\pr\Big( \sum_{i=1}^{n-1} \c_{i,\tbth} ({\bf \Sigma}\otimes {\bf \Sigma}^{-1}) \c_{i,\tbth}\pr\Big)^{-1} \nonumber\\
& \hspace{50mm} \times  \left(\sum_{i=1}^{n-1} \c_{i,\tbth}  (n-i)^{1/2} {\rm vec}({\bGamma_{i; \mathcal{N}}^{(n)}(\bth)})\right).  \label{def.GamStar}
\end{align}

Lemma~\ref{Lemport} below shows that ${\bGamma}_{i, {\cal N}}^{(n)}(\hat{\bth}\n_{\cal N})$ and  ${\bGamma}_{i; {\cal N}}^{(n)*}({\bth})$ under  ${\rm P}^{(n)}_{\tbth ;f}$ only differ by a~$o_{\rm P}(n^{-1/2})$ quantity: under  ${\rm P}^{(n)}_{\tbth ;f}$ and contiguous alternatives, thus,  using~${\bGamma}_{i, {\cal N}}^{(n)}(\hat{\bth}\n_{\cal N})$ in the portmanteau test statistic  is asymptotically equivalent to  using ${\bGamma}_{i; {\cal N}}^{(n)*}({\bth})$. The advantage of the latter is that it no longer involves~$\hat{\bth}\n_{\cal N}$, which simplifies the derivation of asymptotic results. See the Appendix  for a proof.

\begin{lemma}\label{Lemport}
Let Assumptions~\ref{ass.stationary}, \ref{ass.den}, \ref{asylinN}, and \ref{ass.fourMom}  hold. Then 
\begin{equation}\label{equiv*}
(n-i)^{1/2} \text{\rm vec} ({\bGamma}_{i; {\cal N}}^{(n)}(\hat{\bth}\n_{\cal N}) - {\bGamma}_{i; {\cal N}}^{(n)*}({\bth})) = o_{\rm P}(1),
\end{equation}
for any fixed $i\geq 1$, under any  ${\rm P}^{(n)}_{\tbth; f}$  and any contiguous ${\rm P}^{(n)}_{\tbth+ n^{-1/2}\btau; f}$, \linebreak as $i<n\to\infty$.  
\end{lemma}

Note, however, that the asymptotic equivalence \eqref{equiv*} results from the fact that the estimator $\hat\bth_{\cal N}\n$  and the cross-covariances ${\bGamma}_{i; {\cal N}}^{(n)}$ both are associated with the Gaussian distribution (although $f$ needs not to be Gaussian): the same result would not hold with non-Gaussian quasi-maximum likelihood estimators (for instance, the estimator~$\hat\bth_{f}\n$ obtained from the actual likelihood equations $\bDelta\n_{f}(\bth)=\0$, with~$f$ non-Gaussian. 
  Nor does it hold, e.g., with robust or R-estimators (see \cite{HLL2022}) of $\bth$  plugged into~${\bGamma}_{i; {\cal N}}^{(n)}$. 

%
%

\subsection{The (pseudo-)Gaussian test}

The statistic of the Gaussian portmanteau test of \cite{Hosk1980} (to be performed as a pseudo-Gaussian test) takes the form
\begin{align*}
{Q}_{m; {\cal N}}\n(\hat{\bth}\n_{\cal N}) \coloneqq &  \sum_{i=1}^m (n-i) {\rm vec}({\bGamma}_{i, {\cal N}}^{(n)}(\hat{\bth}^{(n)}_{{\cal N}}))\pr (\widehat{\bf \Sigma} \otimes  \widehat{\bf \Sigma})^{-1}  {\rm vec}( {\bGamma}_{i, {\cal N}}^{(n)}(\hat{\bth}^{(n)}_{{\cal N}})) \\
=&  n {\bGamma}_{{\cal N}}^{(m, n)\prime}(\hat{\bth}\n_{\cal N}) (\I_m \otimes \widehat{\bf \Sigma} \otimes  \widehat{\bf \Sigma})^{-1}  {\bGamma}_{{\cal N}}^{(m, n)}(\hat{\bth}\n_{\cal N}),
\end{align*}
where $${\bGamma}_{{\cal N}}^{(m, n)}(\bth) \coloneqq   n^{-1/2} \left( (n-1)^{1/2} \big(\text{vec}( {\bGamma}_{1, {\cal N}}^{(n)}(\bth))\big)^\prime, \ldots,    (n-m)^{1/2} \big(\text{vec}( {\bGamma}_{m, {\cal N}}^{(n)}(\bth))\big)^\prime \right)^\prime$$
and $\widehat{\bf \Sigma} \coloneqq n^{-1} \sum_{t=1}^n \ZZ_t(\hat{\bth}^{(n)}_{{\cal N}}) \ZZ\pr_t(\hat{\bth}^{(n)}_{{\cal N}})$. 

The asymptotic behavior of ${Q}_{m; {\cal N}}\n(\hat{\bth}\n_{\cal N})$ is obtained by  approximating it (for  sufficiently large $m$) by a statistic that is asymptotically chi-square with~$d^2(m-p-q)$ degrees of freedom as $n\to\infty$. Specifically, it follows from Lemma~\ref{Lemport} and the exponential decrease (as $i\to\infty$) of $\Vert\c_{i, \tbth}\Vert$'s that, for  $m$   large enough,   $(n-i)^{1/2} {\rm vec} ({\bGamma}^{(n)}_{i; {\cal N}}(\hat{\bth}\n_{\cal N}))$ is arbitrarily close to 
\begin{align*}
&(n-i)^{1/2} {\rm vec} ({\bGamma}^{(n)**}_{i; {\cal N}}(\bth)) \nonumber \\
&\ = (n-i)^{1/2} {\rm vec} ({\bGamma}^{(n)}_{i; {\cal N}}(\bth)) - ({\bf \Sigma}\otimes {\bf \Sigma}^{-1}) \c_{i,\tbth}\pr\Big( \sum_{i=1}^{m} \c_{i,\tbth} ({\bf \Sigma}\otimes {\bf \Sigma}^{-1}) \c_{i,\tbth}\pr\Big)^{-1} \nonumber\\
& \hspace{50mm} \times  \left(\sum_{i=1}^{m} \c_{i,\tbth}    (n-i)^{1/2} {\rm vec}({\bGamma_{i; \mathcal{N}}^{(n)}(\bth)})\right),
\end{align*}
a quantity that results from truncating at $m$  terms the sums in the right-hand side of~\eqref{def.GamStar}.

Let  
$$
{Q}_{m; {\cal N}}^{(n)**}(\bth) \coloneqq   \sum_{i=1}^m (n-i) ({\rm vec}({\bGamma}_{i, {\cal N}}^{(n)**}(\bth)))\pr ({\bf \Sigma} \otimes  {\bf \Sigma})^{-1}  {\rm vec}( {\bGamma}_{i, {\cal N}}^{(n)**}(\bth)).
$$
The following result  states that ${Q}_{m; {\cal N}}\n(\hat{\bth}^{(n)}_{{\cal N}})$ can be approximated by ${Q}_{m;  {\cal N}}^{(n)**} (\bth)$, where ${Q}_{m;  {\cal N}}^{(n)**} (\bth)$ is asymptotically chi-square under ${\rm P}^{(n)}_{\tbth; f}$ as $n\to\infty$.

\begin{proposition}\label{Prop.GaussAsymptotic}
Let Assumptions~\ref{ass.stationary}, \ref{ass.den}, \ref{asylinN}, and \ref{ass.fourMom} hold. Then, under ${\rm P}^{(n)}_{\tbth; f}$ and contiguous alternatives, 
\begin{enumerate}
\item[{\it (i)}] for all $\delta,\, \varepsilon >0$, there exist integers $M_{\delta, \varepsilon}$, and $N_{\delta, \varepsilon}$ 
  such that 
$${\rm P}\left(\vert {Q}_{m; {\cal N}}\n(\hat{\bth}^{(n)}_{{\cal N}}) - {Q}_{m;  {\cal N}}^{(n)**} (\bth)\vert < \delta \right) > 1 - \varepsilon$$
for all $n \geq N_{\delta, \varepsilon}$ and fixed  $ m \geq M_{\delta, \varepsilon} $;
\item[{\it (ii)}] ${Q}_{m;  {\cal N}}^{(n)**} (\bth)$ is asymptotically chi-square with $d^2(m-p-q)$ degrees of freedom as $n\to\infty$.
\end{enumerate}
\end{proposition}
See the Appendix for a proof.

%

\section{A  portmanteau test based on multivariate center-outward ranks and signs}\label{Sec.RankPortmanteau} 

\subsection{Pseudo-Gaussian versus rank-based tests}

Pseudo-Gaussian tests are generally considered  asymptotically valid over some class $\cal P$ of distributions (containing the Gaussian). And the critical value obtained from the common asymptotic distribution of a pseudo-Gaussian test statistic indeed yields, for each ${\rm P}\in \cal P$, a test with correct asymptotic size~$\alpha$, say. As mentioned in the introduction, the convergence to~$\alpha$ of the finite-sample sizes of the resulting tests, however, in general is highly non-uniform: while  converging to $\alpha $ pointwise in~${\rm P}\in \cal P$, the $\sup_{{\rm P}\in \cal P}$ of these sizes typically fails to do so. If $\rm P$ remains unspecified, which is the case when pseudo-Gaussian tests are performed, the  $\alpha$-level constraint 
 should apply to the $\sup_{{\rm P}\in \cal P}$ of the size under $\rm P$; this $\sup$, however, typically  fails to converge to the asymptotic nominal $\alpha$ level. This is to be kept in mind when performing pseudo-Gaussian tests.

Distribution-free tests do not suffer from the same drawback, since   their finite-sample size   (hence also their asymptotic size)  does not depend on~${\rm P}\in~\!\cal P$, pointwise and  uniform convergence being   equivalent under distribution-free\-ness.~The typical example of distribution-free tests  is that of rank-based tests. The advantages in terms of size and validity of rank tests over the pseudo-Gaussian ones, moreover, are not obtained at the cost of power and efficiency, as shown, for location and simple-output regression, by the celebrated Chernoff-Savage and Hodges-Lehmann results (\cite{ChernoffS, HodgesL}) and, for univariate ARMA time series, by \cite{H94} and \cite{HT00}. Although no fully general versions of the Chernoff-Savage and Hodges-Lehmann inequalities have been established  so far, in the multiple-output case, for center-outward ranks and signs, partial results (restricted to elliptical densities $f$) have been obtained   by \cite{HP02} and \cite{BDS}, and, for elliptical VARMA models, by \cite{HP02',HP05}. These inequalities are quite likely to hold beyond the class of elliptical densities, though.

The basic tools for constructing our rank-based multivariate portmanteau test are the center-outward ranks and signs proposed by \cite{Hallinetal2021}, based on  ideas and results from measure transportation theory. In Section~\ref{secranks}, we introduce the center-outward distribution function and its empirical version, from which the center-outward ranks and signs can be defined. 


\subsection{Center-outward ranks and signs}\label{secranks}

 
Denote by $\mathbb S _d$ and ${\bar{\mathbb S}}_d$ the open and closed unit ball, respectively, and by~${\mathcal S}_{d-1}$ the unit hypersphere   in $\mathbb R ^d$. 
Let $\mathcal{P}_d^\pm$ denote the family of all distributions~$\rm P$ with densities in ${\mathcal F}_d$ such that,  for all positive $r\in\mathbb{R}$, there exist constants~$L^-_r >0$ and $L^+_r <\infty$ for which   
$$L^-_r\leq f({\bf x})\leq L_r^+\quad\text{ for all }\quad {\bf x}\in r\,{\bar{\mathbb S}}_d.$$ 
 For $\rm P$ in  this family,  the center-outward distribution functions defined below  are continuous: see \cite{Figalli2018}. More general cases are studied in \cite{delBarrio2020}, but require  more cautious and less intuitive definitions of these center-outward  functions which, for the sake of simplicity, we do not consider here. Denote by~${\rm U}_d$ the spherical uniform distribution over~${\mathbb S}_d$, that is, the product of a uniform measure over the hypersphere ${\mathcal S}_{d-1}$ and a uniform over the unit interval of distances to the origin. 
 
 The {\it center-outward distribution function}~$\F_{{\pms}}$ of $\rm P$ is defined as the a.e.\ unique gradient of convex function mapping $\R^d$ to $\mathbb{S}_{d}$  and  {\it pushing~$\rm P$ forward} to    ${\rm U}_d$ (that is, such that $\F_{{\pms}}({\bf X})\sim{\rm U}_d$ if ${\bf X}\sim{\rm P}$). For ${\rm P}\in{\mathcal P}_d^\pm$, such mapping is  a homeomorphism between~${\mathbb S}_d\setminus\{{\bf 0}\}$ and~$\mathbb{R}^d\setminus \F_{{\pms}}^{-1}(\{{\bf 0}\})$ (\cite{Figalli2018}) and   the corresponding {\it center-outward quantile function} is defined  as~$\Q_{\pms} \coloneqq \F_{\pms}^{-1}$ (letting, with a small abuse of notation, $\Q_{\pms} ({\bf 0}) \coloneqq \F_{\pms}^{-1}(\{{\bf 0}\})$). For any given distribution~$\rm P$, the quantile fuction $\Q_{{\pms}}$ induces a collection of continuous, connected, and nested quantile  contours $\Q_{\pms}(r{\mathcal S}_{d-1})$ and regions $\Q_{\pms}(r{\mathbb S}_{d})$ of order $r\in[0,1)$;   the {\it center-outward median}  $\Q_{\pms}(\0)$ is a uniquely defined  compact set of Lebesgue measure zero. We refer to \cite{Hallinetal2021} for details.

\begin{figure}[ht!]
\centering
\framebox{{\includegraphics*[angle= -89.5, scale=.350, trim= 0 113 10 110, clip]{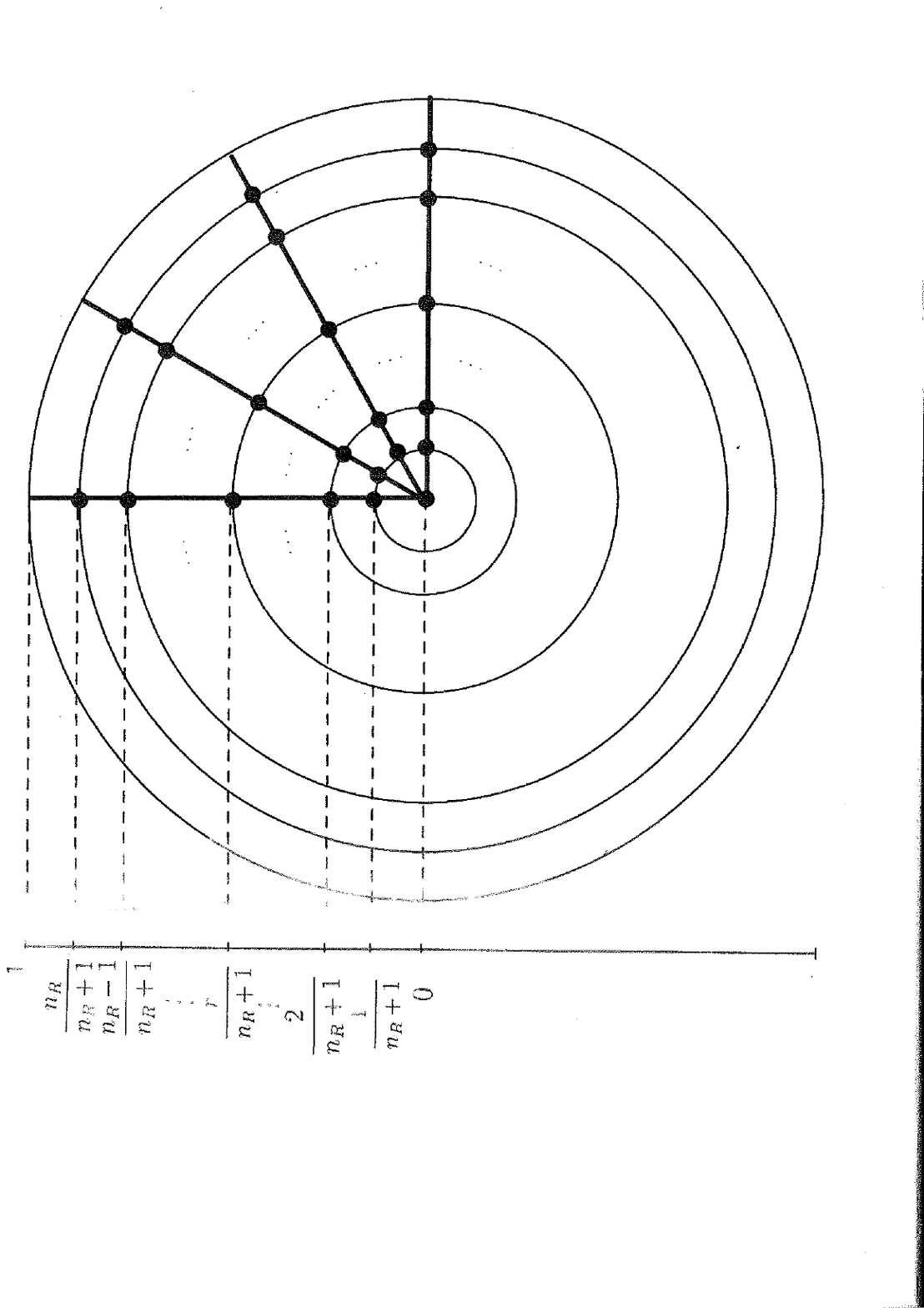}}}
\begin{caption}{\small A regular grid of~$n=n_Rn_S$ points over $\mathbb{S}_2$. }\label{Ffig}
\end{caption}
\end{figure}

 Turning to the sample,  
 the residuals $\ZZ_1^{(n)}(\bth),\ldots,\ZZ_n^{(n)}(\bth)$ under ${\rm P}^{(n)}_{\tbth ;f}$ are i.i.d.\ with density $f\in{\cal F}_d$ and center-outward 
distribution function $\F_\pm$. 
 For the empirical counterpart~$\F^{(n)}_\pm$ of $\F_\pm$,   let $n$ factorize into $n=n_R n_S +n_0,$ for~$n_R, n_S, n_0 \in \mathbb{N}$ and~$0\leq n_0 < \min \{ n_R, n_S \}$, where~$n_R \rightarrow \infty$ and $n_S \rightarrow \infty$ as~$n \rightarrow \infty$, and consider a sequence $\mathfrak{G}\n$ 
  of grids, where each grid
consists of the $n_Rn_S$ intersections between an $n_S$-tuple $(\boldsymbol{u}_1,\ldots \boldsymbol{u}_{n_S})$ of unit vectors, and the~$n_R$ hyperspheres  with radii $1/(n_R+1),\ldots ,n_R/(n_R+1)$ centered at the origin, along with~$n_0$ copies of the origin. The only requirement is that 
 the discrete distribution with probability masses~$1/n$ at each gridpoint and probability mass~$n_0/n$ at the origin converges weakly to the uniform~${\rm U}_d$ over the ball $\mathbb{S}_d$. 
 Then, we  define~$\F_{\pm}^{(n)}(\ZZ^{(n)}_t) $, for~$t = 1,\ldots ,n$ as the solution (optimal mapping) of a coupling problem between the residuals 
and the gridpoints. Specifically, the empirical center-outward distribution function is the (random) discrete mapping\vspace{-2mm} 
$$
\F^{(n)}_{\pm}: \ZZ^{(n)}\coloneqq (\ZZ^{(n)}_1,\ldots ,\ZZ^{(n)}_n) \mapsto (\F^{(n)}_{\pm}(\ZZ^{(n)}_1),\ldots ,\F^{(n)}_{\pm}(\ZZ^{(n)}_n))
$$
satisfying 
\begin{equation}\label{Fpm0}
\s \Vert \ZZ^{(n)}_t - \F^{(n)}_\pm(\ZZ^{(n)}_t)\Vert ^2 = \underset{T \in \mathcal{T}}{\min} \s \Vert \ZZ^{(n)}_t - T(\ZZ^{(n)}_t)\Vert ^2,
\end{equation}
where $\ZZ^{(n)}_t=\ZZ^{(n)}_t(\bth)$, the set $\{\F^{(n)}_{\pm}(\ZZ^{(n)}_t)\vert t = 1,\ldots ,n\}$ coincides with the $n$ points of the grid, and $\mathcal{T}$ denotes the set of all possible bijective mappings between $\ZZ^{(n)}_1, \ldots , \ZZ^{(n)}_n$ and the $n$ gridpoints. Intuition for this mapping in dimension $d=2$ is provided in Figure~\ref{Ffig}.

Based on this empirical center-outward distribution function,  
  the {\it center-outward ranks} are defined as 
\begin{equation}
{R^{(n)}_{\pm, t}}\coloneqq  {R^{(n)}_{\pm, t}}(\bth)\coloneqq({n_R + 1}) \Vert \F^{(n)}_\pm(\ZZ^{(n)}_t) \Vert,  \label{Ranks}
\end{equation}
 the {\it center-outward signs}  as 
\begin{equation}
\S^{(n)}_{\pm, t}\coloneqq  {\S^{(n)}_{\pm, t}}(\bth)\coloneqq  \F^{(n)}_\pm(\ZZ^{(n)}_t)  I [\F^{(n)}_\pm(\ZZ^{(n)}_t) \neq \boldsymbol{0}]/ 
\Vert \F^{(n)}_\pm(\ZZ^{(n)}_t) \Vert. \label{Signs}
\end{equation}
 It follows that $\F^{(n)}_\pm(\ZZ_t^{(n)})$ factorizes into\vspace{-2mm}
\begin{equation}\label{factFpm}
\F^{(n)}_\pm(\ZZ_t^{(n)})=   \frac{R^{(n)}_{\pm, t}}{n_R + 1}  \S^{(n)}_{\pm, t},\quad\text{ whence }\quad \ZZ_t^{(n)}= \Q^{(n)}_\pm\Big(
 \frac{R^{(n)}_{\pm, t}}{n_R + 1}  \S^{(n)}_{\pm, t}\Big).
\end{equation}

Those ranks and signs are jointly distribution-free  under~${\rm P}^{(n)}_{\tbth; f}$ (for  ${\rm P}\!\in~\!{\cal P}_d^\pm$): more precisely,  under~${\rm P}^{(n)}_{\tbth; f}$,  the $n$-tuple $\F^{(n)}_\pm(\ZZ_1^{(n)}),\ldots , \F^{(n)}_\pm(\ZZ_n^{(n)})$ is uniformly distributed over the $n!/n_0!$ permutations with repetition of the~$n$ underlying  gridpoints (the origin having multiplicity~$n_0$). Moreover, the center-outward distribution functions, ranks, and signs inherit, from the invariance properties of Euclidean distances, elementary but remarkable invariance and equi\-variance properties with respect to shift, global scale, and orthogonal transformations: see Proposition~2.2 in   \cite{HHH22}   for details.

An intuitive choice of $n_R$ and $n_S$ in the factorization of $n$  is  $n_R\approx n^{1/d}$ and~$n_S\approx n^{(d-1)/d}$;  $n_R$, indeed,  is the cardinality of a one-dimensional grid over~$[0,1]$,  $n_S$   the cardinality of a grid over the $(d-1)$-sphere ${\mathcal S}_{d-1}$. Other heuristic criteria are possible, though. \cite{GM21},  for instance, since the grid is supposed to provide an approximation of the spherical uniform, suggests   minimizing the Wasserstein distance between the empirical distribution over the grid  and the spherical uniform: see Section~3.2 of \cite{HM2022}.  

\subsection{Center-outward rank-based cross-covariance matrices} \label{SecSR}



Writing $\F^{(n)}_{\pm, t}$, $R^{(n)}_{\pm, t}$, and $\S^{(n)}_{\pm, t}$   for~$\F^{(n)}_\pm(\ZZ_t^{(n)}(\bth))$, $ {R^{(n)}_{\pm, t}}(\bth)$, and~$ {\S^{(n)}_{\pm, t}}(\bth)$, respectively, consider the center-outward rank-based counterpart of ${\bGamma}^{(n)}_f (\bth)$. Specifically, define \vspace{-2mm}  
 \begin{align}
\tenq{\bGamma}_{\J_1, \J_2}^{(n)}(\bth) \coloneqq &   n^{-1/2} \left( (n-1)^{1/2} \big(\text{vec}(  \tenq{\bGamma}_{1, \J_1, \J_2}^{(n)}(\bth))\big)^\prime, \ldots , \right. \nonumber \\
&\ \  \left.  (n-i)^{1/2} \big(\text{vec}(  \tenq{\bGamma}_{i, \J_1, \J_2}^{(n)}(\bth))\big)^\prime, \ldots ,
\big(\text{vec}( \tenq{\bGamma}_{n-1, \J_1, \J_2}^{(n)}(\bth))\big)^\prime\right)^\prime\!, \label{tilde.S.m}
\end{align}
\vspace{-10mm}

\noindent with, for $ i = 1, \ldots , n - 1,$  \vspace{-2mm} 
 \begin{equation}\label{tildeGam1}
\tenq{\bGamma}_{i, \J_1, \J_2}^{(n)}(\bth) \coloneqq (n-i)^{-1} \sum_{t=i+1}^n \J_1\left(\frac{R^{(n)}_{\pm, t}}{n_R + 1}\S^{(n)}_{\pm, t}\right) \J_2\pr\left(\frac{R^{(n)}_{\pm, t-i}}{n_R + 1}\S^{(n)}_{\pm, t-i}\right)   
\end{equation}
where  $\J_1$ and~$ \J_2: {\mathcal S}_d \rightarrow \R$
 are score functions satisfying Assumption~\ref{ass.vp} below. Call $\tenq{\bGamma}_{i, \J_1, \J_2}^{(n)}(\bth) $  a (residual)  {\it  rank-based~cross-covariance matrix with lag $i$}: this matrix  is distribution-free under ${\rm P}\n_{\tbth; f}$ due to distribution-freeness of the center-outward ranks and signs.  In order to establish its asymptotic distribution, 
   we make the following assumption on $\J_1$ and~$ \J_2$.

\begin{assumption}\label{ass.vp}
 The score functions {\rm $\J_1$ and $\J_2$
 \begin{compactenum}
\item[{\it (i)}]  are continuous over $\mathbb{S}_d$;
\item[{\it (ii)}]   are  square-integrable, that is, 
$
 \int_{\mathbb{S}_d} \Vert  \J_\ell (\u) \Vert^2 {\rm d U}_d < \infty$ for $ \ell =1, 2$.
 \end{compactenum} 
Moreover,
  \begin{compactenum}
\item[{\it (iii)}]   for any sequence $\mathfrak{s}^{(n)} \coloneqq \{{\bf s}_1^{(n)}, \ldots, {\bf s}_n^{(n)}\}$ of $n$-tuples in  
$\mathbb{S}_d$ such that the uniform discrete distribution 
over $\mathfrak{s}^{(n)}\!$ converges weakly to ${\rm U}_d$, 
\begin{equation}\label{ass2}
\underset{n\rightarrow\infty}{\lim} n^{-1} \sum_{t=1}^n 
 \Vert \J_\ell ({\bf s}_t^{(n)}) \Vert^2 =    \int_{\mathbb{S}_d} \Vert  \J_\ell (\u) \Vert^2 {\rm d U}_d  , \quad \ell =1, 2.
\end{equation}
\end{compactenum}}
\end{assumption}

 Let 
 \begin{align*}{\mub}\n &\coloneqq {\rm E}(\J_1(\F\n_{\pm, 2}) \J\pr_2(\F\n_{\pm, 1})),\\ 
{\D}_{\J_1, \J_2}&\coloneqq   \left[ \int_{\mathbb{S}_d} \J_2(\u_2) \J_2\pr(\u_2) {\rm dU}_d(\u_2) \right] \otimes \left[ \int_{\mathbb{S}_d} \J_1(\u_2) \J_1\pr(\u_2) {\rm dU}_d(\u_2) \right]  \\ 
& - \left[ \int_{\mathbb{S}_d} \J_2(\u) {\rm dU}_d   \int_{\mathbb{S}_d} \J_2\pr(\u) {\rm dU}_d \right] \otimes \left[  \int_{\mathbb{S}_d} \J_1(\u) {\rm dU}_d   \int_{\mathbb{S}_d} \J_1\pr(\u) {\rm dU}_d  \right],
 \end{align*}
and 
\begin{align*}
\K_{\J_1, \J_2, f} &\coloneqq   {\rm E}[\J_2(\F_{\pms}(\bepsilon_t)) \bepsilon_t\pr] \otimes {\rm E}[\J_1(\F_{\pms}(\bepsilon_t)) \bvp^\prime_{f}(\bepsilon_t)].
\end{align*}
We then have the following result from~\cite{HLL2022}.

 \begin{proposition}\label{asy.Gami.utilde}
Let Assumptions~\ref{ass.stationary}, \ref{ass.den}, and~\ref{ass.vp} hold. Then, 
 for any posi\-tive~integers $i\neq j$, the vectors
 $$(n-i)^{1/2}\text{\text{\rm vec}} (\tenq{\bGamma}_{i, {\J}_1, {\J}_2}^{(n)}(\bth) - \mub^{(n)})\quad\text{and}\quad (n-j)^{1/2}\text{\text{\rm vec}} (\tenq{\bGamma}_{j, {\J}_1, {\J}_2}^{(n)}(\bth) - \mub^{(n)})$$
  are jointly asymptotically normal, 
with mean $(\0\pr, \0\pr)\pr$ under  ${\rm P}^{(n)}_{\tbth; f}$,  mean\vspace{-1mm}
\begin{equation}\label{aslingam}
\left((\K_{\J_1, \J_2, f} \c_{i, \tbth}\pr \btau)\pr, (\K_{\J_1, \J_2, f} \c_{j, \tbth}\pr \btau)\pr\right)\pr
\end{equation}
under ${\rm P}^{(n)}_{\tbth + n^{-1/2}\btau ;f}$, and   covariance  $\left(
\begin{array}{cc} {\D}_{\J_1, \J_2}& \0 \\ \0& {\D}_{\J_1, \J_2} \end{array}
\right)$  under both.
\end{proposition}

Without loss of generality, we will assume that $\mub\n = o(n^{-1/2})$. A sufficient condition is either $\int_{\mathbb{S}_d} \J_1(\u){\rm d}U_d = {\boldsymbol 0}$ or $\int_{\mathbb{S}_d} \J_2(\u){\rm d}U_d =  {\boldsymbol 0}$; indeed,  following the same lines as in Lemma~1 of \cite{HL2020},  one has  
$$\mub\n - {\rm E}(\J_1(\F_{\pm, 1})) {\rm E}(\J\pr_2(\F_{\pm, 1})) = o(n^{-1/2}).$$
  Moreover, note that,  for the scores in Example 1-3 of Section~\ref{Sec.Example}, $\mub\n = \0$ holds whenever the regular grid $\mathfrak{G}\n$ is symmetric with respect to the origin. See \cite{HLL2023} for  details.

To define the rank-based test statistic and derive its asymptotic theory, we assume the following asymptotic linearity of $(n-i)^{1/2} {\rm vec} (\tenq{\bGamma}^{(n)}_{i, \J_1, \J_2}(\bth))$, where the form of the right-hand side of \eqref{assA4} again follows from the form of the shift matrix in \eqref{aslingam}.

\begin{assumption}\label{asylin1}
For any positive integer $i$, as $n\to\infty$, \vspace{-2mm}
\begin{equation}\label{assA4}
(n-i)^{1/2}\left[{\rm vec} (\tenq{\bGamma}^{(n)}_{i, \J_1, \J_2}(\bth + n^{-1/2}\btau)) - (\tenq{\bGamma}^{(n)}_{i, \J_1, \J_2}(\bth)) \right] = - \K_{\J_1, \J_2, f} \c_{i, \tbth}\pr \btau + o_{\rm P}(1)
\end{equation}
\vspace{-6mm}

\noindent under  ${\rm P}^{(n)}_{\tbth ;f}$ (hence also under ${\rm P}^{(n)}_{\tbth + n^{-1/2}\btau ;f}$).
\end{assumption}
Sufficient conditions for \eqref{assA4}  can be found in \cite{vdAetal2015}. 

Let $\tenq{\bDelta}^{(n)}_{\J_1, \J_2}(\bth) \coloneqq \sum_{i=1}^{n-1} \c_{i,\tbth}  (n-i)^{1/2} {\rm vec}(\tenq{\bGamma}^{(n)}_{i, \J_1, \J_2}(\bth))$ denote the {\it rank-based central sequence}. The asymptotic linearity of 
$\tenq{\bDelta}^{(n)}_{\J_1, \J_2}(\bth)$,  
 that is,
\begin{equation}\label{asylin2}
\tenq{\bDelta}^{(n)}_{\J_1, \J_2}(\bth+ n^{-1/2}\btau) - \tenq{\bDelta}^{(n)}_{\J_1, \J_2}(\bth) = -  \sum_{i=1}^{n-1}  \c_{i, \tbth} \K_{\J_1, \J_2, f} \c\pr_{i, \tbth}\btau + o_{\rm P}(1).
\end{equation}
as $n\to\infty$ readily follows from Assumption~\ref{asylin1}, both under  ${\rm P}^{(n)}_{\tbth ;f}$ and under~${\rm P}^{(n)}_{\tbth + n^{-1/2}\btau ;f}$.


\subsection{Center-outward rank-based  portmanteau  tests}\label{subsec: ranktest}

Recall from Lemma~\ref{Lemport} that the pseudo-Gaussian portmanteau test relies on a  plug-in  of the QMLE $\hat{\bth}\n_{{\cal N}}$, which kind of neutralizes the asymptotic impact of computing  cross-covariances from estimated residuals rather than exact ones.  Lemma~\ref{Lemport.rank} below shows that a similar result holds for the rank-based cross-covariance matrix $\tenq{\bGamma}_{i; {\bf J}_1,{\bf J}_2}^{(n)}$ when   computed at the R-estimator $\tenq{\bth}^{(n)}_{\J_1, \J_2}$ of \cite{HLL2022}. That result will not hold if the R-estimator and the  rank-based cross-covariance matrix are based on distinct score  functions, or if the QMLE is plugged in instead of the adequate R-estimator.

Let ${\bUpsilon}\n_{\J_1, \J_2, f}(\bth) \coloneqq \sum_{i=1}^{n-1}  \c_{i, \tbth} \K_{\J_1, \J_2, f} \c^{\prime}_{i, \tbth}$
denote the shift matrix in \eqref{asylin2} and let $\widehat{\bUpsilon}\n_{\J_1, \J_2}$ be a consistent (under ${\rm P}^{(n)}_{\tbth ;f}$) estimator thereof; see \cite{HLL2022} for details.  
Further, denote by  $\bar{\bth}\n$ a preliminary $\sqrt{n}$-consistent and asymptotically discrete  estimator of $\bth$.  Recall that an estimator $\bar{\bth}\n$ of~$\bth$ is called  {\it asymptotically discrete} if the number of values it   takes in balls of radius  $cn^{-1/2}$ ($c>0$) is bounded as $n\to\infty$, which is easily obtained by discretization. This assumption, which is classical in Le Cam asymptotics, is needed in asymptotic statements but has no finite-sample implications. The {\it one-step R-estimator} of $\bth$ is defined as
\begin{equation}\label{def.Rest}
\tenq{\bth}^{(n)}_{\J_1, \J_2} \coloneqq \bar{\bth}\n + n^{-1/2} (\widehat{\bUpsilon}\n_{\J_1, \J_2})^{-1} \tenq{\bDelta}^{(n)}_{\J_1, \J_2}(\bar{\bth}\n).
\end{equation}

{Now, consider  the matrix $\tenq{\bGamma}_{i; {\bf J}_1, {\bf J}_2, f}^{(n)*}(\bth)$ of residuals in the regression \linebreak  of~$(n-i)^{1/2}\text{\rm vec}\tenq{\bGamma}_{i; {\bf J}_1{\bf J}_2}^{(n)}(\bth)$ with respect to $\tenq{\bDelta}_{{\bf J}_1, {\bf J}_2}^{(n)}(\bth)$ in  the  covariance matrix (note, however, that this matrix is not the asymptotic covariance matrix\linebreak  of~$\left((n-i)^{1/2}\text{\rm vec}\tenq{\bGamma}_{i; {\bf J}_1{\bf J}_2}^{(n)}(\bth),\,\tenq{\bDelta}_{{\bf J}_1, {\bf J}_2}^{(n)}(\bth)\right)$)
\[
\left(
\begin{array}{cc}
{\bf D}_{\J_1, \J_2}&\K_{\J_1, \J_2, f} \c_{i, \tbth}^{\prime}\\
\c_{i, \tbth}\K_{\J_1, \J_2, f} &\displaystyle \left(\sum_{i=1}^{n-1} \c_{i, \tbth} \K_{\J_1, \J_2, f} \c^{\prime}_{i, \tbth}\right)
\end{array}
\right).
\]
  Namely, define $(n-i)^{1/2}\tenq{\bGamma}_{i; {\bf J}_1, {\bf J}_2, f}^{(n)*}(\bth)$ as }
\begin{align}
&(n-i)^{1/2} {\rm vec} (\tenq{\bGamma}^{(n)*}_{i, \J_1, \J_2, f}(\bth)) \nonumber \\
&\coloneqq (n-i)^{1/2} {\rm vec} (\tenq{\bGamma}^{(n)}_{i, \J_1, \J_2}(\bth)) - \K_{\J_1, \J_2, f} \c_{i, \tbth}^{\prime}  \left(\sum_{i=1}^{n-1} \c_{i, \tbth} \K_{\J_1, \J_2, f} \c^{\prime}_{i, \tbth}\right)^{-1} \tenq{\bDelta}^{(n)}_{\J_1, \J_2}(\bth) \nonumber \\
&= (n-i)^{1/2} {\rm vec} (\tenq{\bGamma}^{(n)}_{i, \J_1, \J_2}(\bth)) - \K_{\J_1, \J_2, f} \c_{i, \tbth}^{\prime}  (\sum_{i=1}^{n-1} \c_{i, \tbth} \K_{\J_1, \J_2, f} \c^{\prime}_{i, \tbth})^{-1} \nonumber\\
& \hspace{50mm} \times (\sum_{i=1}^{n-1} \c_{i, \tbth}  (n-i)^{1/2} {\rm vec}(\tenq{\bGamma}_{i, \J_1, \J_2}\n(\bth))).  \label{def.GamJ1J2Star}
\end{align}
Lemma~\ref{Lemport.rank} below shows that $(n-i)^{1/2}\tenq{\bGamma}_{i, \J_1, \J_2}^{(n)}(\tenq{\bth}_{\J_1, \J_2}^{(n)})$ and the oracle statistic~$(n-i)^{1/2}\tenq{\bGamma}_{i, \J_1, \J_2, f}^{(n)*}(\tenq{\bth}_{\J_1, \J_2}^{(n)})$ are asymptotically equivalent. Hence, the asymptotics of the rank-based test statistic  constructed from $(n-i)^{1/2}\tenq{\bGamma}_{i, \J_1, \J_2}^{(n)}(\tenq{\bth}_{\J_1, \J_2}^{(n)})$ coincide with those of  its counterpart built from $(n-i)^{1/2}\tenq{\bGamma}_{i, \J_1, \J_2, f}^{(n)*}(\tenq{\bth}_{\J_1, \J_2}^{(n)})$. See the Appendix for a proof.

\begin{lemma}\label{Lemport.rank}
Let Assumptions~\ref{ass.stationary}, \ref{ass.den}, \ref{ass.vp}, and~\ref{asylin1} hold. Then
\begin{enumerate}
\item[{\it (i)}]$\tenq{\bDelta}^{(n)}_{\J_1, \J_2}(\tenq{\bth}^{(n)}_{\J_1, \J_2})=o_{\rm P}(1)$,
\item[{\it (ii)}]   $(n-i)^{1/2} {\rm vec} (\tenq{\bGamma}_{i, \J_1, \J_2, f}^{(n)*}(\tenq{\bth}^{(n)}_{\J_1, \J_2}) - \tenq{\bGamma}_{i, \J_1, \J_2}^{(n)}(\tenq{\bth}^{(n)}_{\J_1, \J_2})) = o_{\rm P}(1)$, and 
\item[{\it (iii)}]   $(n-i)^{1/2} {\rm vec} (\tenq{\bGamma}_{i, \J_1, \J_2, f}^{(n)*}(\tenq{\bth}^{(n)}_{\J_1, \J_2}) - \tenq{\bGamma}_{i, \J_1, \J_2, f}^{(n)*}({\bth})) = o_{\rm P}(1),$
\end{enumerate}
for any fixed $i\geq 1$, under  ${\rm P}^{(n)}_{\tbth; f}$  and ${\rm P}^{(n)}_{\tbth + n^{1/2}\btau; f}$, as $n\to\infty$.
\end{lemma}
\vspace{3mm}

Note that the asymptotic covariance matrix of $(n-i)^{1/2}\text{\text{\rm vec}} (\tenq{\bGamma}_{i, {\J}_1, {\J}_2}^{(n)}(\bth))$ in Proposition~\ref{asy.Gami.utilde}, which does not depend on the underlying density $f$, differs from the corresponding shift matrix~$\K_{\J_1, \J_2, f}$, which is not even symmetric for a general $f$. This was not the case for  
 the asymptotic covariance matrix of $(n-i)^{1/2}\text{\text{\rm vec}} ({\bGamma}_{i, {\cal N}}^{(n)}(\bth))$. As a consequence, the construction of rank-based portmanteau tests is more intricate than in the classical pseudo-Gaussian case and requires some additional notation. 

In view of Lemma~\ref{Lemport.rank}, the test statistic should involve a consistent estimator of the covariance matrix of  $(n-i)^{1/2}{\rm vec} ({\bGamma}^{(n)*}_{i; {\J}_1, {\J}_2}({\bth}))$, rather than a consistent estimator of $(n-i)^{1/2}{\rm vec} (\tenq{\bGamma}^{(n)}_{i; {\J}_1, {\J}_2}({\bth}))$. Let 
\begin{equation}\label{Omegadef}{\bOmega}^{(m)}_{i, \J_1, \J_2, f}(\bth) \coloneqq {\W}^{(m)}_{i, \J_1, \J_2, f}(\bth) {\W}^{(m)\prime}_{i, \J_1, \J_2, f}(\bth)
\end{equation}
 with 
\begin{align}\label{Wdef}
{\W}^{(m)}_{i, \J_1, \J_2, f}(\bth) &\coloneqq (\e^{(m)\prime}_i \otimes \D_{\J_1, \J_2}^{1/2}) - \K_{\J_1, \J_2, f} \c_{i, \tbth}^{\prime}  (\sum_{i=1}^{m} \c_{i, \tbth} \K_{\J_1, \J_2, f} \c^{\prime}_{i, \tbth})^{-1} \\
& \hspace{60mm} \times \C^{(m+1)}_{\tbth} (\I_{m} \otimes \D_{\J_1, \J_2}^{1/2})\nonumber
\end{align}
and  ${\e}_i^{(m)}$  standing for the $i$th vector of the canonical basis in~${\mathbb{R}}^{m}$.

For   sufficiently large $m$, due to exponential decrease of $\Vert\c_{i, \tbth}\Vert$, the matrices~${\bOmega}^{(m)}_{i, {\J}_1, {\J}_2, f}(\tenq{\bth}^{(n)}_{\J_1, \J_2})$ are (uniformly in $i$)  arbitrarily close to the covariance matrices of  $(n-i)^{1/2}{\rm vec} ({\bGamma}^{(n)*}_{i; {\J}_1, {\J}_2, f}({\bth}))$ as $n\rightarrow \infty$.
Letting
\begin{equation}\label{Edef}
\E^{(m)}_{\J_1, \J_2, f}(\bth) \coloneqq \I_{md^2} - (\I_m \otimes  \K_{\J_1, \J_2, f}) \C^{(m+1)\prime}_{\tbth} (\sum_{i=1}^{m} \c_{i, \tbth} \K_{\J_1, \J_2, f} \c^{\prime}_{i, \tbth})^{-1} \C^{(m+1)}_{\tbth},
\end{equation}
we have 
\begin{equation}\label{eq.diagOmega}
{\rm diag}({\bOmega}^{(m)}_{i, \J_1, \J_2, f}(\bth))_{1\leq i \leq m} = \E^{(m)}_{\J_1, \J_2, f}(\bth) (\I_m \otimes \D_{\J_1, \J_2}) \E^{(m)\prime}_{\J_1, \J_2, f}(\bth).
\end{equation}
Heuristically, $\E^{(m)}_{\J_1, \J_2, f}(\bth)$, which appears in the rank-based test statistic \eqref{eq:tildeQ2}, is taking into account the difference between $\tenq{\bGamma}^{(n)*}_{i, \J_1, \J_2, f}(\bth)$ and $\tenq{\bGamma}^{(n)}_{i, \J_1, \J_2}(\bth)$ in~\eqref{def.GamJ1J2Star}.


Our rank-based test statistics involve  consistent estimators of the quantities defined above. Denote  by $\widehat{\K}_{\J_1, \J_2}\n$ a consistent estimator  of $\K_{\J_1, \J_2, f}$; such~$\widehat{\K}\n_{\J_1, \J_2}$ can be obtained via~\eqref{assA4}---for example, letting
$$\btau_j = -\c_{1, \tenq{\tbth}^{(n)}_{\J_1, \J_2}} (\c_{1, \tenq{\tbth}^{(n)}_{\J_1, \J_2}}\pr \c_{1, \tenq{\tbth}^{(n)}_{\J_1, \J_2}})^{-1}  {\bf e}_j^{(d^2)},$$ where ${\bf e}_j^{(d^2)}$ denotes $j$th vector of the canonical basis in $\mathbb{R}^{d^2}$, then
$$(n-1)^{1/2}\left[{\rm vec} ({\bGamma}^{(n)}_{1; \J_1, \J_2}(\tenq{\bth}^{(n)}_{\J_1, \J_2} + n^{-1/2}\btau_j)  - {\bGamma}^{(n)}_{1; \J_1, \J_2}(\tenq{\bth}^{(n)}_{\J_1, \J_2})) \right]$$ is a consistent estimator of the $j$th column of $\K_{\J_1, \J_2, f}$. For $j = 1, \ldots,~d^2$, this yields a consistent estimator of  $\K_{\J_1, \J_2, f}$. Let~$\widehat{\bOmega}^{(m)}_{i, \J_1, \J_2}(\bth)$ and $\widehat{\E}^{(m)}_{\J_1, \J_2}(\bth)$ stand for the consistent estimators   of  ${\bOmega}_{i, \J_1, \J_2, f}(\bth)$ and~$\E^{(m)}_{\J_1, \J_2, f}(\bth)$ obtained by plugging $\widehat{\K}_{\J_1, \J_2}$ into \eqref{Omegadef}--\eqref{Wdef} and \eqref{Edef}. Then the center-outward rank-based portmanteau test statistic  we are proposing takes the form
\begin{align}
\tenq{Q}_{m; \J_1, \J_2}\n(\tenq{\bth}^{(n)}_{\J_1, \J_2}) &\coloneqq 
\sum_{i=1}^m (n-i) {\rm vec}(\tenq{\bGamma}_{i, \J_1, \J_2}^{(n)}(\tenq{\bth}^{(n)}_{\J_1, \J_2}))\pr (\widehat{\bOmega}^{(m)}_{i, \J_1, \J_2}(\tenq{\bth}^{(n)}_{\J_1, \J_2}))^{-}\nonumber  \\ 
&\hspace{50mm} \times {\rm vec}( \tenq{\bGamma}_{i, \J_1, \J_2}^{(n)}(\tenq{\bth}^{(n)}_{\J_1, \J_2})), 
\end{align}
where $\A^{-}$ denotes the Moore-Penrose inverse of a matrix $\A$. Alternatively, in view of \eqref{eq.diagOmega}, $\tenq{Q}_{m; \J_1, \J_2}\n(\tenq{\bth}^{(n)}_{\J_1, \J_2}) $ can be written as
\begin{align}\label{eq:tildeQ2}
\tenq{Q}_{m; \J_1, \J_2}\n(\tenq{\bth}^{(n)}_{\J_1, \J_2}) 
&= n \tenq{\bGamma}_{\J_1, \J_2}^{(m, n)\prime}(\tenq{\bth}^{(n)}_{\J_1, \J_2}) \left(\widehat{\E}^{(m)}_{\J_1, \J_2}(\tenq{\bth}^{(n)}_{\J_1, \J_2}) (\I_m \otimes \D_{\J_1, \J_2}) \widehat{\E}^{(m)\prime}_{\J_1, \J_2}(\tenq{\bth}^{(n)}_{\J_1, \J_2})\right)^{-} \nonumber \\ 
&\hspace{55mm} \times \tenq{\bGamma}_{\J_1, \J_2}^{(m, n)}(\tenq{\bth}^{(n)}_{\J_1, \J_2})
\end{align}\vspace{-12mm}

\noindent where
\begin{align*}\tenq{\bGamma}_{\J_1, \J_2}^{(m, n)}(\bth) \coloneqq   n^{-1/2} \left( (n-1)^{1/2} \right.&(\text{vec}( \tenq{\bGamma}_{1, \J_1, \J_2}^{(n)}(\bth)))^\prime,  \ldots \\
&\left. \qquad\qquad\ldots,   (n-m)^{1/2} (\text{vec}(  \tenq{\bGamma}_{m, \J_1, \J_2}^{(n)}(\bth)))^\prime \right)^\prime
\end{align*}
is a truncated version of $\tenq{\bGamma}_{\J_1, \J_2}^{(n)}(\bth)$. 

\subsection{Asymptotic distribution}

Due to Lemma~\ref{Lemport.rank} and the consistency of $\tenq{\bth}^{(n)}_{\J_1, \J_2}$ and $\widehat{\K}_{\J_1, \J_2}$, we have 
\begin{equation}\label{eq.PPstar}
\tenq{Q}_{m; \J_1, \J_2}\n(\tenq{\bth}^{(n)}_{\J_1, \J_2}) = \tenq{Q}_{m; \J_1, \J_2, f}^{(n)*} (\bth) + o_{\rm P}(1),
\end{equation}
where
$$\tenq{Q}_{m; \J_1, \J_2, f}^{(n)*} (\bth) \coloneqq \sum_{i=1}^m (n-i) {\rm vec}(\tenq{\bGamma}_{i, \J_1, \J_2, f}^{(n)*}(\bth))\pr ({\bOmega}^{(m)}_{i, \J_1, \J_2; f}(\bth))^{-}  {\rm vec}(\tenq{\bGamma}_{i, \J_1, \J_2, f}^{(n)*}(\bth)).$$
Moreover, note that for $m$ large enough, due again to the exponential decrease of $\Vert\c_{i, \tbth}\Vert$, $\tenq{Q}_{m; \J_1, \J_2, f}^{(n)*} (\bth)$ is arbitrarily close to 
\begin{align}
\tenq{Q}_{m; \J_1, \J_2, f}^{(n)**} (\bth) &\coloneqq \sum_{i=1}^m (n-i) {\rm vec}(\tenq{\bGamma}_{i, \J_1, \J_2, f}^{(m, n)**}(\bth))\pr ({\bOmega}^{(m)}_{i, \J_1, \J_2; f}(\bth))^{-} 
{\rm vec}(\tenq{\bGamma}_{i, \J_1, \J_2, f}^{(m, n)**}(\bth)) \label{Pstarstar}
\end{align}
\vspace{-12mm}

\noindent where 
\begin{align}
&(n-i)^{1/2} {\rm vec} (\tenq{\bGamma}^{(m, n)**}_{i, \J_1, \J_2, f}(\bth)) \label{Gamstarstar} \\ 
&\hspace{20mm}\coloneqq (n-i)^{1/2} {\rm vec} (\tenq{\bGamma}^{(n)}_{i, \J_1, \J_2}(\bth)) - \K_{\J_1, \J_2, f} \c_{i, \tbth}^{\prime}  (\sum_{i=1}^{m} \c_{i, \tbth} \K_{\J_1, \J_2, f} \c^{\prime}_{i, \tbth})^{-1} \nonumber \\
& \hspace{68mm} \times (\sum_{i=1}^{m} \c_{i, \tbth}  (n-i)^{1/2} {\rm vec}(\tenq{\bGamma}_{i, \J_1, \J_2}\n(\bth)))\nonumber\
\end{align}
\vspace{-7mm}

\noindent is an approximation of $(n-i)^{1/2} {\rm vec} (\tenq{\bGamma}^{(n)*}_{i, \J_1, \J_2, f}(\bth))$ resulting from truncating at~$m$ the summation in the right hand side of~\eqref{def.GamJ1J2Star};    the asymptotic  covariance matrix of $(n-i)^{1/2} {\rm vec} (\tenq{\bGamma}^{(m, n)**}_{i, \J_1, \J_2, f}(\bth))$ is ${\bOmega}^{(m)}_{i, \J_1, \J_2; f}(\bth)$ under ${\rm P}\n_{\tbth; f}$.

The following result  is the rank-based counterpart of Proposition~\ref{Prop.GaussAsymptotic} and states that $\tenq{Q}_{m; \J_1, \J_2}\n(\tenq{\bth}^{(n)}_{\J_1, \J_2})$, for $m$ large enough, can be approximated by~$\tenq{Q}_{m; \J_1, \J_2, f}^{(n)**} (\bth)$ which is asymptotically chi-square under ${\rm P}\n_{\tbth;f}$ as $n\to\infty$.  
Note that $\tenq{\bGamma}_{\J_1, \J_2}^{(m, n)}(\bth)$, which depends on the ranks of the ``true'' residuals, is fully distribution-free, while (due to the presence of $ \K_{\J_1, \J_2, f} $) $\tenq{\bGamma}^{(m, n)**}_{i, \J_1, \J_2, f}(\bth)$ in \eqref{Gamstarstar}, hence 
 $\tenq{Q}_{m; \J_1, \J_2, f}^{(n)**} (\bth)$ in \eqref{Pstarstar}, are not. However, $\tenq{Q}_{m; \J_1, \J_2, f}^{(n)**} (\bth)$ is asymptotically distribution-free. 

%

\begin{proposition}\label{Prop.ranktest}
Let Assumptions~\ref{ass.stationary}, \ref{ass.den}, \ref{ass.vp}, and~\ref{asylin1} hold. Then, under ${\rm P}^{(n)}_{\tbth; f}$, 
\begin{enumerate}
\item[{\it (i)}] for all $\delta,\, \varepsilon >0$, there exist $M_{\delta, \varepsilon}$ and $N_{\delta, \varepsilon} \in \mathbb{N}$ such that 
$${\rm P}\left(\vert \tenq{Q}_{m; \J_1, \J_2}\n(\tenq{\bth}^{(n)}_{\J_1, \J_2}) - \tenq{Q}_{m; \J_1, \J_2, f}^{(n)**} (\bth) \vert < \delta \right) > 1 - \varepsilon$$
for all $n \geq N_{\delta, \varepsilon}$ and $M_{\delta, \varepsilon} \leq m \leq n-1$;
\item[{\it (ii)}] irrespective of the innovation density $f$, $\tenq{Q}_{m; \J_1, \J_2, f}^{(n)**} (\bth)$ is asymptotically chi-square with $d^2(m-p-q)$ degrees of freedom as $n\to\infty$.
\end{enumerate}
\end{proposition}
See the Appendix for a proof.

\subsection{Local powers and AREs}
Having obtained the asymptotic distributions of portmanteau test statistics  under the null,  natural questions are: what are their local asymptotic powers? can we derive their asymptotic relative efficiencies (AREs) with respect, for instance, to the classical procedure? 

Local powers and AREs, however, depend on the type of alternative under consideration and are usually obtained via an application of Le Cam's third Lemma. This requires an embedding of the VARMA$(p,q)$ model that has been studied so far into a larger model enjoying LAN and containing the alternatives of interest.  A simple example is the VARMA$(p_1,q_1)$ model with~$p_1>p$ and/or $q_1>q$, for which LAN results are available \citep{GH95}; call it the ``VARMA$(p_1,q_1)$ case''. Other alternatives of interest are, e.g.,  the presence of bilinear terms in the data-generating process, or the conditional heteroskedasticity of the ${\boldsymbol\epsilon}_t$'s---for which LAN, in the vicinity of VARMA$(p,q)$ models, is not available in the literature. While we have no doubt that such LAN structures hold  under appropriate assumptions,  establishing such results is beyond the scope of this paper, and are restricting this section to a discussion of the VARMA$(p_1,q_1)$ case. 

Now, even in the VARMA$(p_1,q_1)$ case,  things are not as simple as it may appear at first sight. While asymptotic shifts under local VARMA$(p_1,q_1)$ alternatives are easily obtained 
for the Gaussian $(n-i)^{1/2}$vec$\left(\tenq{\bGamma}_{i, \J_1, \J_2, f}^{(n)*}({\bth})\right)$ and the rank-based~$(n-i)^{1/2}$vec$\left(\tenq{\bGamma}_{i, \J_1, \J_2, f}^{(n)*}({\bth})\right)$, hence, in view of Lemmas~\ref{Lemport} and~\ref{Lemport.rank}, 
for $(n-i)^{1/2}$vec$\left({\bGamma}_{i; {\cal N}}^{(n)}(\hat{\bth}\n_{\cal N})\right)$ and~$(n-i)^{1/2}$vec$\left(\tenq{\bGamma}_{i, \J_1, \J_2}^{(n)}(\tenq{\bth}^{(n)}_{\J_1, \J_2})\right)$, the quadratic forms ${Q}_{m; {\cal N}}^{(n)*}({\bth})$ and $\tenq{Q}_{m; \J_1, \J_2}^{(n)*}({\bth})$ are not {\it idempotent}   (not even asymptotically so) in these statistics and are not asymptotically chi-square under the null.   As a consequence, the asymptotic shifts of $(n-i)^{1/2}$vec$\left({\bGamma}_{i; {\cal N}}^{(n)}(\hat{\bth}\n_{\cal N})\right)$ and~$(n-i)^{1/2}$vec$\left(\tenq{\bGamma}_{i, \J_1, \J_2, f}^{(n)*}(\tenq{\bth}^{(n)}_{\J_1, \J_2})\right)$ do not induce, for the portmanteau sta\-tis\-tics~${Q}_{m; {\cal N}}\n(\hat{\bth}^{(n)})\!=\!{Q}_{m; {\cal N}}^{(n)*}({\bth})+~\!o_{\rm P}(1)$ and $\tenq{Q}_{m; \J_1, \J_2}\n(\tenq{\bth}^{(n)}_{\J_1, \J_2})\!=\! \tenq{Q}_{m; \J_1, \J_2}^{(n)*}({\bth})+~\!o_{\rm P}(1)$, the usual 
chi-square noncentrality parameters, and the classical methods expressing AREs in terms of ratios of noncentrality parameters do not apply. 

These classical methods, on the other hand, do apply to the asymptotically chi-square approximations ${Q}_{m; {\cal N}}^{(n)**}(\bth)$ and $\tenq{Q}_{m; \J_1, \J_2}^{(n)**}({\bth})$. These, however, are oracle test statistics that cannot be computed from the observations and the resulting AREs  cannot   be considered as approximations of the actual portmanteau AREs; the approximation errors in Propositions~\ref{Prop.GaussAsymptotic} and~\ref{Prop.ranktest}, indeed, do not go to zero as $n\to\infty$. Such AREs are of limited practical value,  thus. Accordingly, we will not proceed any further with their numerical evaluation, and rather recommend, for power comparisons, the Monte Carlo approach adopted in Section~\ref{Sec:Numerical}. 

\subsection{Some standard score functions}\label{Sec.Example}

The rank-based test statistic $\tenq{Q}_{m; \J_1, \J_2}\n(\tenq{\bth}^{(n)}_{\J_1, \J_2})$ depends on the choice of the score functions $\J_1$ and $\J_2$. Here we provide three examples of standard score functions, extending scores that are widely applied in the univariate  (see e.g. \cite{HL2020}) and  the elliptical   setting  (see \cite{HP04}). 

\textbf{Example 1} ({\it Sign test} scores). Setting 
$\J_{\ell}\left(\frac{R^{(n)}_{\pm, t}}{n_R + 1}\S^{(n)}_{\pm, t}\right)  = \S^{(n)}_{\pm, t}, \quad \ell = 1, 2$  
yields the center-outward sign-based  cross-covariance matrices\vspace{-2mm} 
\begin{align*}
\tenq{\bGamma}_{i, \text{sign}}^{(n)}(\bth) = (n-i)^{-1} \sum_{t=i+1}^n   \S^{(n)}_{{\pms}, t}(\bth) \S^{(n)\prime}_{{\pms}, t-i}(\bth), \quad i = 1, \ldots , n-1. 
\end{align*}
\vspace{-6mm}

\noindent The resulting test statistic $\tenq{Q}_{m; \J_1, \J_2}\n(\tenq{\bth}^{(n)}_{\J_1, \J_2})$  entirely relies on the center-outward  signs $\S^{(n)}_{{\pms}, t}(\bth)$, which justifies the terminology  {\it sign test}  scores.\medskip

\textbf{Example 2} ({\it Spearman} scores).   Another simple choice is $\J_1(\u)  = \u =~\!\J_2(\u)$.  
The corresponding rank-based cross-covariance matrices  are\vspace{-2mm}
\begin{align*}
\tenq{\bGamma}_{i, \text{Sp}}^{(n)}(\bth) = (n-i)^{-1} \sum_{t=i+1}^n   \F^{(n)}_{{\pms}, t} \F^{(n)\prime}_{{\pms}, t-i}, \quad i = 1, \ldots , n-1
\end{align*}
\vspace{-4mm}

\noindent reducing, for $d=1$,   to Spearman autocorrelations, whence the terminology  {\it Spearman}~scores. \medskip


 \textbf{Example 3} ({\it van der Waerden} or spherical  {\it normal} scores). Let
$$\J_{\ell} \left(\frac{R^{(n)}_{\pm, t}}{n_R + 1}\S^{(n)}_{\pm, t}\right) = J_{\ell} \left(\frac{R^{(n)}_{\pm, t}}{n_R + 1}\right) \S^{(n)}_{\pm, t}, \quad  \ell  = 1, 2,
$$
with $J_{\ell}(u) = \big((F^{\chi^2}_d)^{-1} (u)\big)^{1/2}\!$, where $F^{\chi^2}_d$ denotes the chi-square distribution function with~$d$ degrees of freedom. This yields  the {\it spherical van der Waerden} (vdW) {\it rank scores}, with cross-covariance matrices\vspace{-2mm}  
\begin{align*}
&\tenq{\bGamma}_{i, \text{vdW}}^{(n)}(\bth) = (n-i)^{-1}\!\! \\ 
\hspace{12mm} &  \qquad\times \sum_{t=i+1}^n\! \left[\big(F^{\chi^2}_d\big)^{-1}\! \left(\frac{R^{(n)}_{{\pms}, t}(\bth)}{n_R + 1}\right)\right]^{1/2}\! \left[\big(F^{\chi^2}_d\big)^{-1}\! \left(\frac{R^{(n)}_{{\pms}, t-i}(\bth)}{n_R + 1}\right)\right]^{1/2} 
\!\! \!\! \nonumber \\
\hspace{12mm} &  \qquad\qquad\qquad\qquad\qquad\qquad\times \S^{(n)}_{{\pms}, t} (\bth)\S^{(n)\prime}_{{\pms}, t-i}(\bth),  
\qquad i = 1, \ldots , n-1.  
\end{align*}

\section{Numerical assessment of finite-sample performance and resistance to outliers}\label{Sec:Numerical}

In this section, we investigate through a brief Monte Carlo experiment the finite-sample performance of the center-outward rank-based and pseudo-Gaus\-sian~tests under various innovation densities, including contaminated~ones.

\subsection{Size and power}\label{Sec:sizepower}
  We first compare the sizes of these tests, hence their validity,  by gene\-rating $N=300$ replications of sample size $n=1000$ from the bivariate VARMA($1, 1$) model with 
 \begin{equation}\label{VARMA(1,1)}
 \A_1 = \left(\begin{array}{cc} 0.5 &0.2 \\ -0.1&0.4
  \end{array}
  \right)
\quad\text{ and } 
\B_1 =\left(\begin{array}{cc} 0.3&0 \\ 0&0.4
  \end{array}
  \right)
\end{equation}
and 
 three types of innovation densities $f$: spherical normal,  mixture of three Gaussians,   and skew-$t$ distribution with $3$ degrees of freedom (denoted as skew-$t_3$; see~\cite{AC03} for a definition), respectively.

\begin{figure}[htbp]\vspace{-15mm}
\begin{subfigure}
\centering 
\includegraphics[width=130mm, height=65mm]{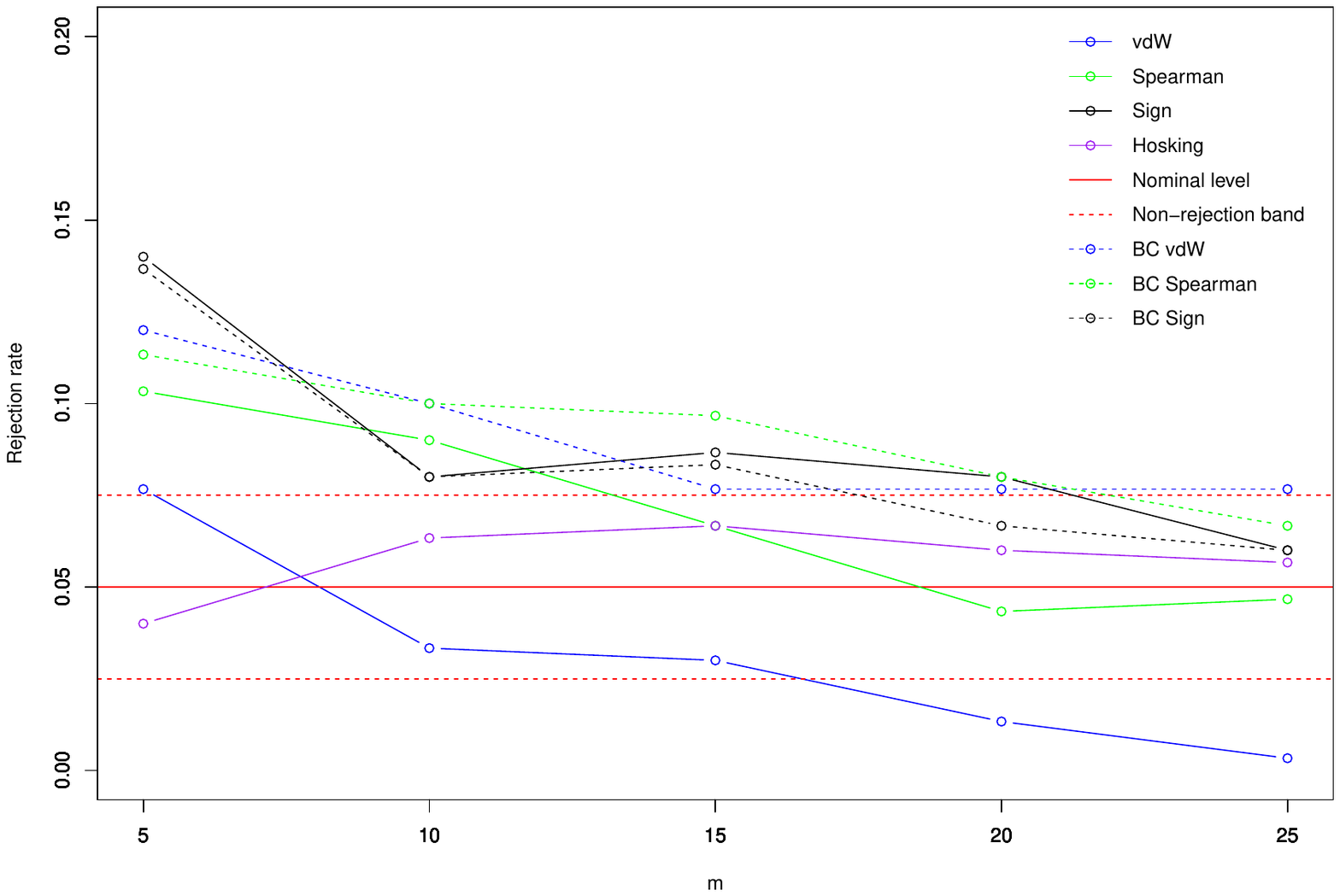}\vspace{-5mm}
\end{subfigure}
\begin{subfigure}
\centering
\includegraphics[width=130mm, height=65mm]{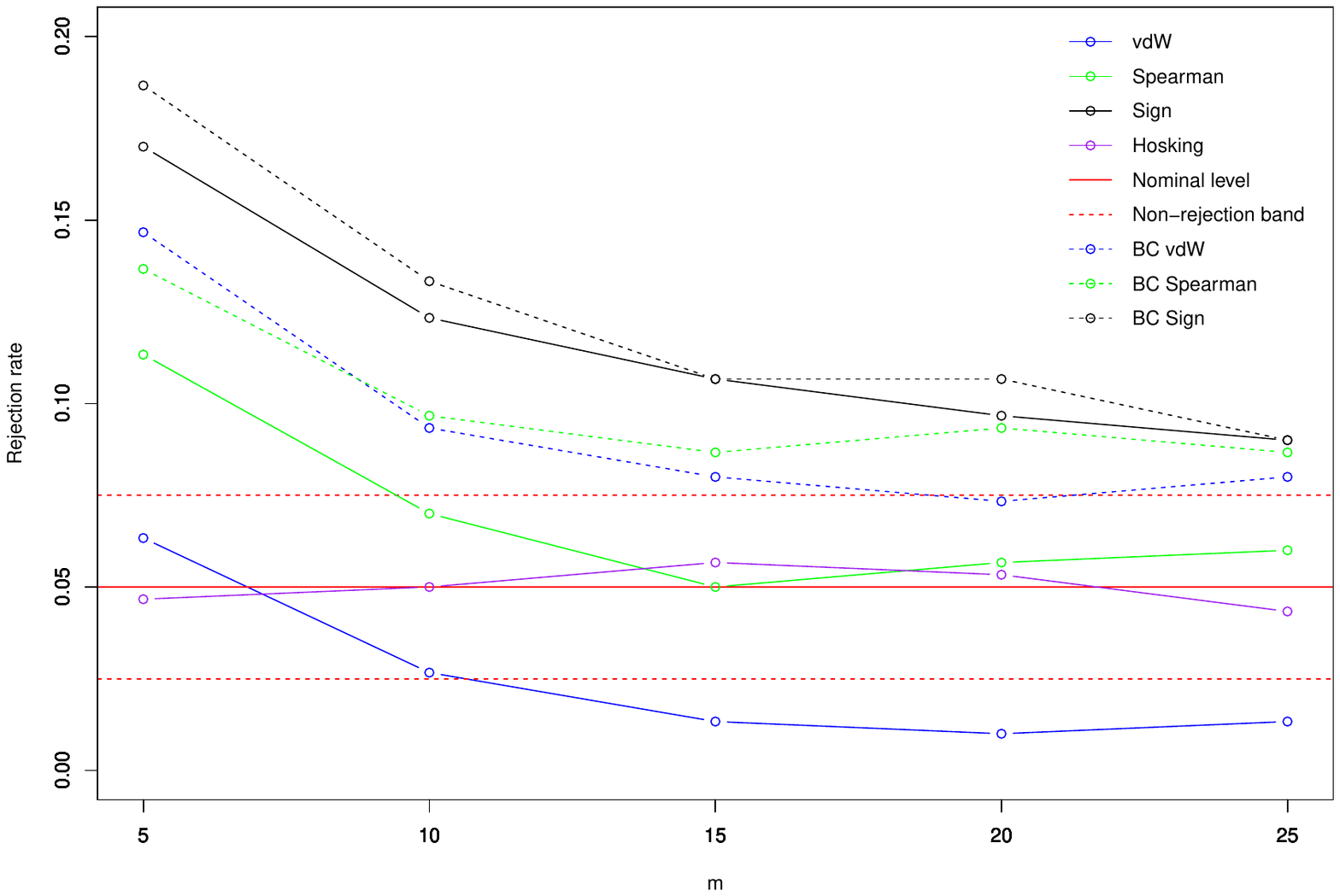}\vspace{-5mm}
\end{subfigure}
\begin{subfigure}
\centering
\includegraphics[width=130mm, height=65mm]{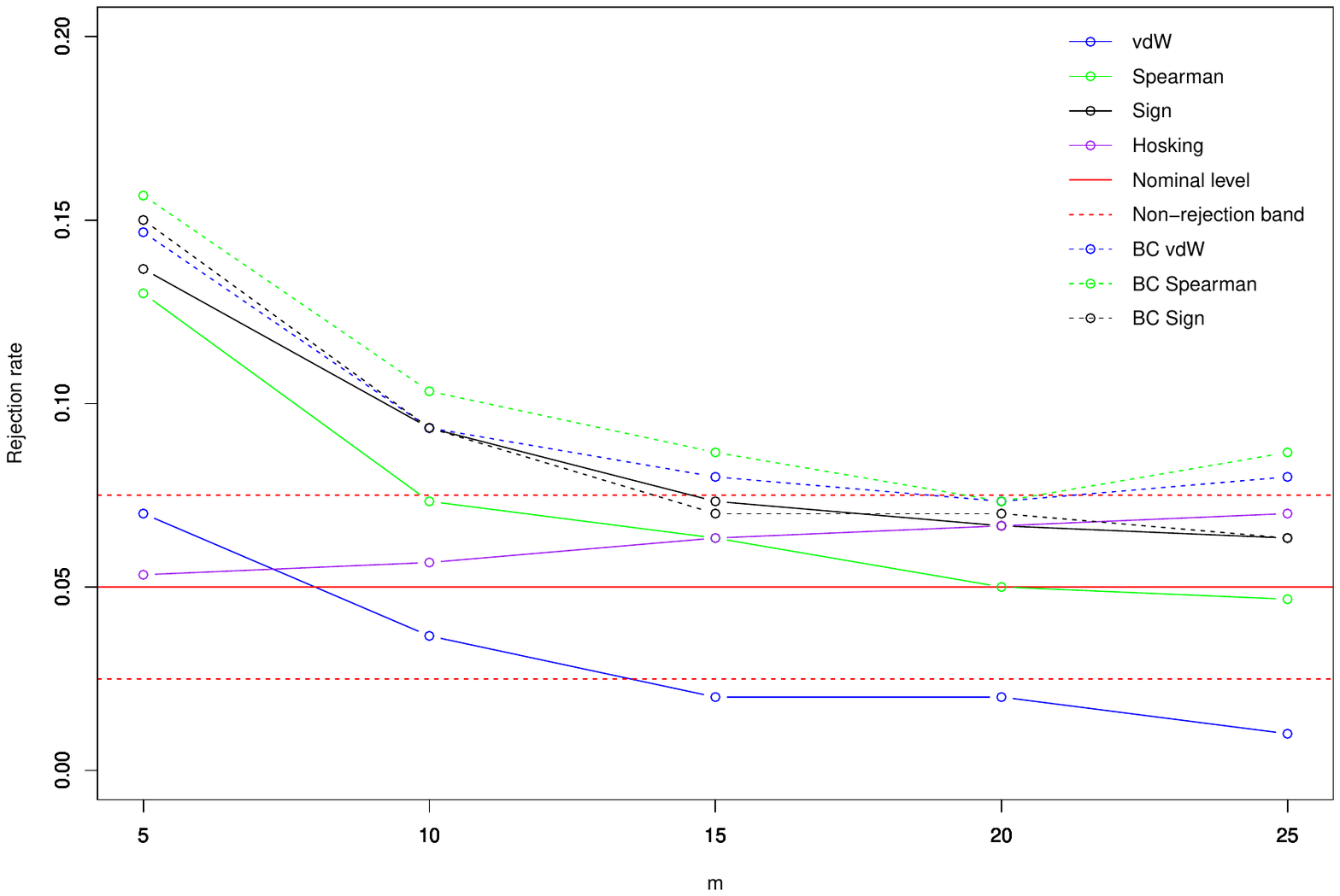}\vspace{-8mm}
\end{subfigure}
\caption{\small Rejection frequencies (nominal level 5\%; asymptotic chi-square critical values and, for    the van der Waerden, Spearman, and sign tests, bias-corrected (BC) permutational critical values),  for~$m =~\!5,\, 10, \ldots, \,25$, 
 of the  Gaussian, van der Waerden, Spearman, and sign portmanteau tests for unspecified VARMA(1,1) model, under the VARMA(1,1) model \eqref{VARMA(1,1)} with  spherical normal (upper panel),  mixture  \eqref{Eq. Mixture} of three Gaussians (middle panel) and  skew-$t_3$  (lower panel) innovation densities.  Number of observations $n=1000$;  $N=300$ repli\-cations. The solid and dashed horizontal lines indicate the nominal level~$\alpha=5\%$ and the rejection limits of the~$5\%$ two-sided test  of the hypothesis that the actual level indeed is~$5\%$. 
}\label{plot.Size}
%
%
%
\end{figure}

 \medskip

 The mixture is of the form\vspace{-2mm} 
\begin{equation}
\frac{3}{8}   {\cal N}(\bmu_1, \bSigma_1) + \frac{3}{8}   {\cal N}(\bmu_2, \bSigma_2) + \frac{1}{4} {\cal N}(\bmu_3, \bSigma_3),
\label{Eq. Mixture}
\end{equation}
with $\bmu_1 = (-5, 0)^\prime,\ \bmu_2 = (5, 0)^\prime, \ \bmu_3 = (0, 0)^\prime$,  and\vspace{-1mm}  
$$\bSigma_1 = 
\begin{pmatrix}
7 & 5\vspace{-1mm} \\
5 & 5\vspace{-1mm}
\end{pmatrix}, \  
\bSigma_2 = 
\begin{pmatrix}
7 &  -6\vspace{-1mm} \\
6 &  6\vspace{-1mm}
\end{pmatrix}, \ 
\bSigma_3 = 
\begin{pmatrix}
4 & 0\vspace{-1mm} \\
0 & 3\vspace{-1mm}
\end{pmatrix}.
\vspace{-2mm}$$

For each replication, a VARMA(1,1) model has been estimated, via Gaussian likelihood, van der Waerden, Spearman,  and sign R-estimation, respectively (ranks and signs computed from a grid with $n_R = 25, n_S = 40$, and~$n_0 =~\! 0$). 
      Based on these estimators, the Gaussian, van der Waerden,  Spearman, and sign portmanteau tests (ranks and signs computed from a grid with $n_R = 25, n_S = 40$, and~$n_0 =~\! 0$) were performed for each replication at $5\%$ nominal level for $m = 5, 10, \ldots, 25$. Note that a series length of $n=1000$ in dimension two---roughly corresponding to a series length of~$n\approx \sqrt{1000}\approx32$ in dimension one---is not a very large one. 
 
The rejection frequencies under the null hypothesis of an unspeci- \linebreak fied~VARMA(1,1) model  are shown in Figure~\ref{plot.Size} for the spherical Gaussian,  the mixture of three Gaussians \eqref{Eq. Mixture}, and  the skew-$t_3$ innovation densities. The dotted horizontal lines provide the $5\%$ critical band outside which the empirical size of a test is significantly different from the nominal level $\alpha =5\%$. 

Irrespective of the  innovation density,   rejection rates (for fixed $n$)  tend to decrease as $m$ increases.  As a function of $m$,   the decrease depends on the scores and the actual innovation. It appears, for instance, that the van der Waerden test meets the nominal $\alpha$-level constraint for $m\approx 8$ under all three innovation densities, while the Spearman test requires $m$ to be as large as 10 (under the mixture of Gaussians or the Skew-$t_3$) or even~14 (under the spherical Gaussian) for not exhibiting significant over-rejection. In that respect, van der Waerden  is ``more $m$-parsimonious'' than Spearman, and  the sign test is particularly ``greedy'' (never meeting the nominal $5\%$ size requirement for $m$ between 5 and 25 in the mixture case).   As for Hosking's traditional pseudo-Gaussian test, it is meeting the nominal  level constraint for $m$ as small as~$5$.

\begin{figure}[btp]\vspace{-13mm}
\begin{subfigure}
\centering
\includegraphics[width=130mm, height=65mm]{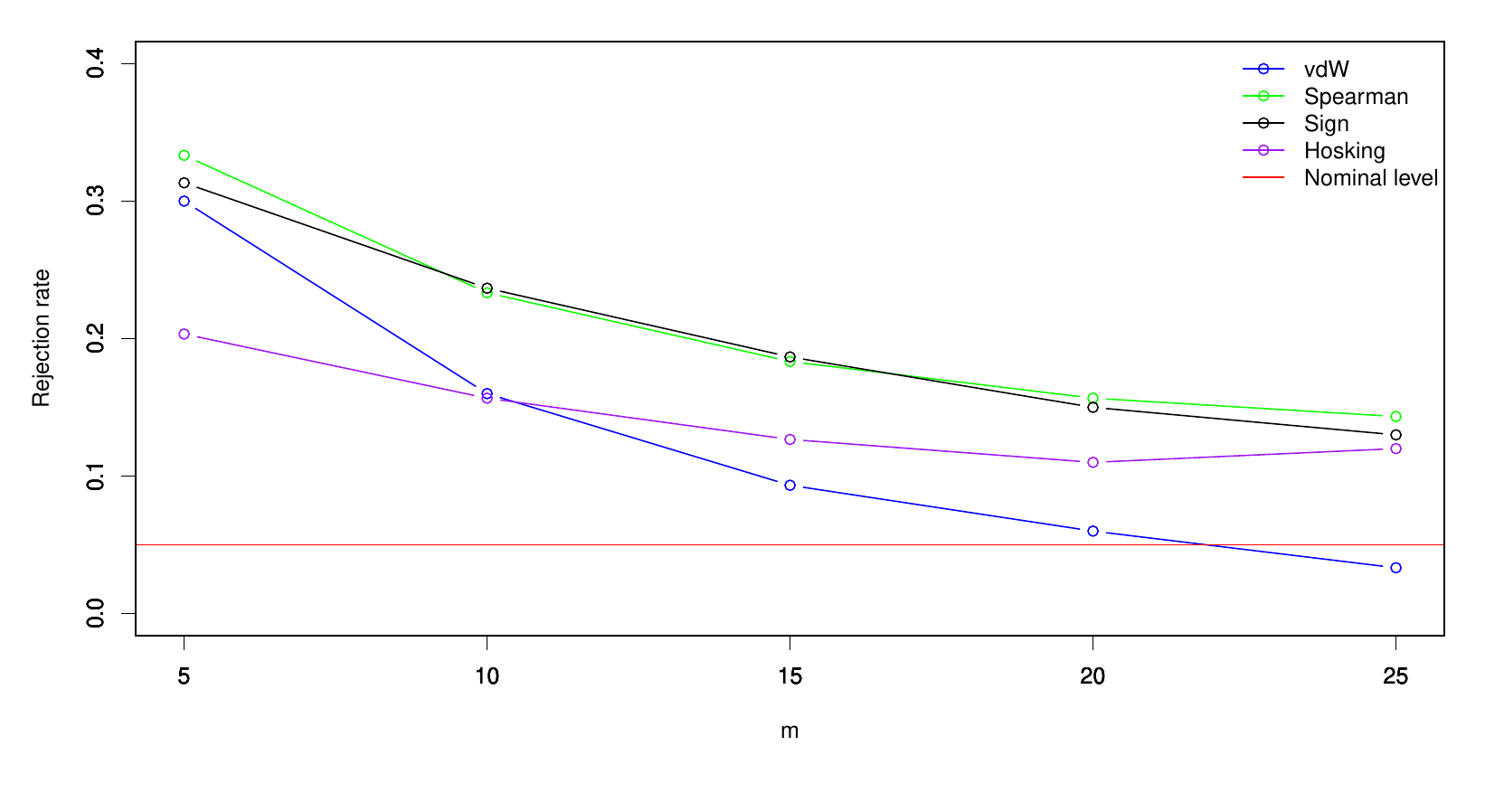}\vspace{-5mm}
\end{subfigure}
\begin{subfigure}
\centering
\includegraphics[width=130mm, height=65mm]{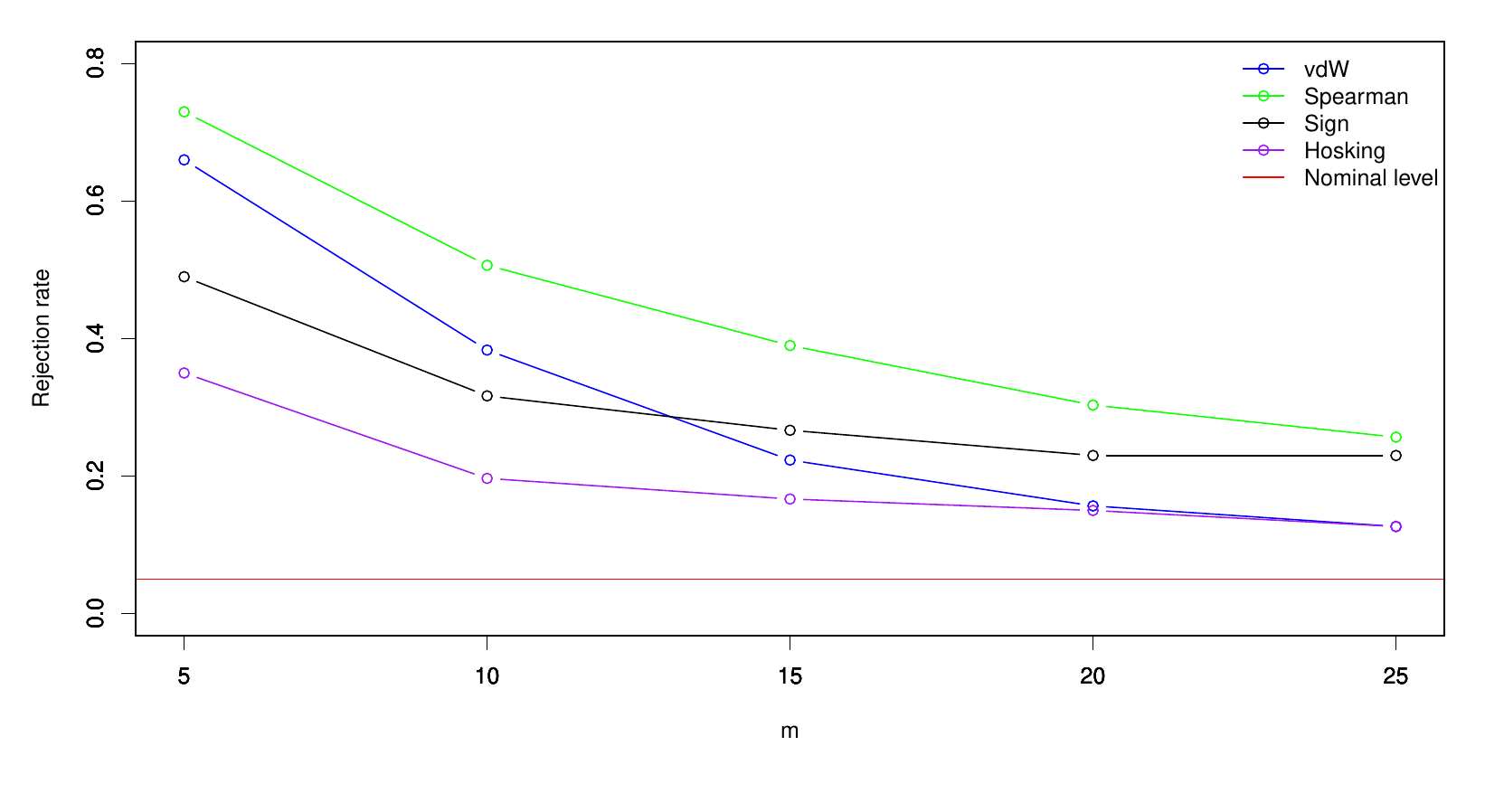}\vspace{-5mm}
\end{subfigure}
\begin{subfigure}
\centering
\includegraphics[width=130mm, height=65mm]{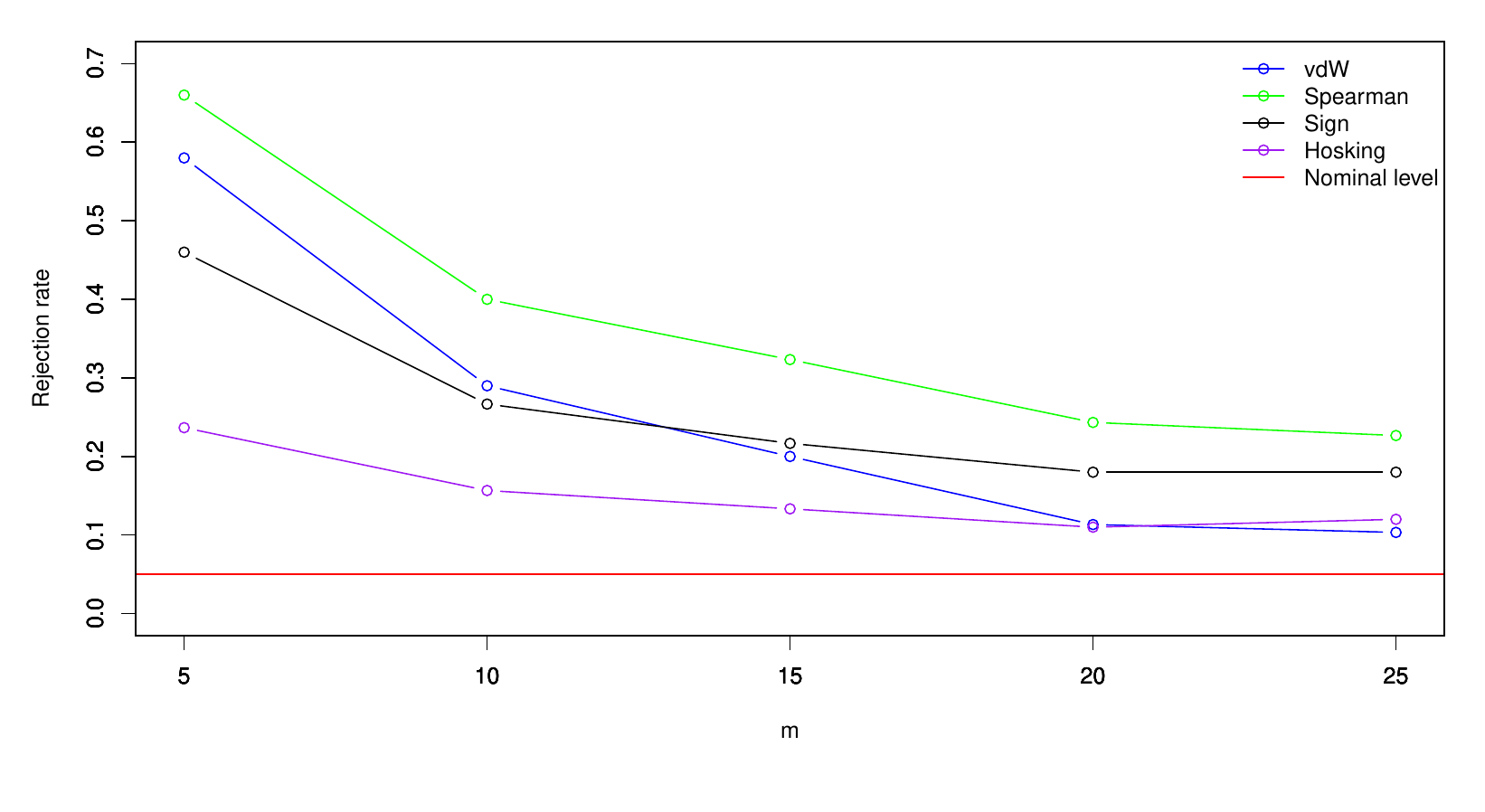}\vspace{-10mm}
\end{subfigure}

\caption{\small Rejection frequencies (nominal level 5\%; asymptotic chi-square critical values),   for~$m\!=~\!5,\,10, \ldots,\, 25$, 
 of the  Gaussian, van der Waerden, Spearman, and sign portmanteau tests for unspecified VARMA(1,1) model, under the  VARMA(1,2) alternative \eqref{alternative} with spherical normal (upper panel),  mixture  \eqref{Eq. Mixture} of three Gaussians (middle panel) and  skew~$t_3$  (lower panel) innovation densities. Number of observations~$n=1000$;  $N=300$ repli\-cations. The solid  hori\-zontal line indicates the nominal level~$\alpha=5\%$. 
}\label{plot.Power}
\end{figure}

In order to compare the powers of the various portmanteau  tests,   we simulated the bivariate VARMA($1, 2$) model with matrix coefficients 
\begin{equation}\label{alternative}
\A_1  \!=\!  \left(\!\begin{array}{cc} 0.5 &0.2 \\ -0.1&0.4
  \end{array}
 \! \right),\quad\!\!\!\!
\B_1 \!=\! \left(\!\begin{array}{cc} 0.3&0 \\ 0&0.4
  \end{array}
 \! \right), \quad\!\!\!\! \text{and}\quad\!\! \!\!\B_2  \!=\!  \left(\!\begin{array}{cc} 0.07&0.03 \\  -0.02&0.1
  \end{array}
\!  \right)
\end{equation}
 under the same innovation densities as above.

The rejection  frequencies, under  the VARMA(1,2) alternative \eqref{alternative} with spherical normal,   mixture  \eqref{Eq. Mixture} of three Gaussians, and   skew-$t_3$ innovation densities, of the same  Gaussian, van der Waerden,  Spearman, and sign portmanteau tests as in Figure~\ref{plot.Size} are shown in Figure~\ref{plot.Power}  for $m = 5, 10, \ldots, 25$. As expected,  rejection frequencies tend to decrease as $m$ increases,  irrespective of the innovation density. Interestingly, for given $m$,  Spearman  has greater power than Hosking's traditional test, even under  normal innovations,  for all~$m \leq25$.  So has the sign test---but this is largely due to the fact that, as shown in  Figure~\ref{plot.Size}, it  is over-rejecting. Now, a fair comparison should take the results of Figure~\ref{plot.Size} into account. Table~\ref{table1}  provides, for each test and each innovation density, the optimal lag numbers $m_0$---that is, the lag compatible with the $5\%$ level condition  at which the empirical power is maximal---along with that   power.  Inspection of the table reveals that van der Waerden with~$m=5$ or $6$ lags is the uniform winner, and more parcimonious than Spearman, which is second best, albeit with more than $m=10$ lags.  Both are overperforming Hosking's traditional test.

%
%
%
%
%
%
%
%
%
%
%
%
\begin{table}[!ht]

\centering

\begin{tabular}{cccccc ccccccc}

\hline\hline

             &    &      \multicolumn{3}{c}{Gaussian} & &    \multicolumn{3}{c}{Mixture}             &  &     \multicolumn{3}{c}{Skew-$t_3$} \\ \cline{3-5} \cline{7-9} \cline{11-13} 

              Test & & $m_0$ & & Power & & $m_0$ & & Power & & $m_0$ &  & Power \\ \cline{1-1} \cline{3-3} \cline{5-5} \cline{7-7} \cline{9-9}  \cline{11-11} \cline{13-13} 

              vdW & & {$6$} & & {${\bf 0.272}$} & & $5$ & & ${\bf 0.660}$ & & $5$ & & ${\bf 0.580}$  \\

              Spearman & & $14$ & & $0.193$ & & $10$ & & $0.507$ & & $10$ & & $0.400$  \\

              Sign & & {$22$} & & $0.142$ & & -- & & -- & & $15$ & & $0.217$  \\

              Hosking & & $5$ & & $0.203$ & & $5$ & & $0.350$ & & $5$ & & $0.237$  \\ \hline \hline

\end{tabular}

\caption{\small Optimal lag numbers $m_0$ and empirical powers (against the VARMA alternative~\eqref{alternative}) of the van der Waerden, Spearman, sign and Hosking portmanteau tests, under spherical normal, mixture of three Gaussians, and skew~$t_3$ innovation densities, respectively. Boldface indicates the winner (consistently, van der Waerden with $m=5$ or $6$) in each column.\label{table1}}\vspace{-3mm}
\end{table}

These conclusions, of course, are based on a very limited Monte Carlo experiment. More general innovation densities and more general alternatives (bilinear, heteroskedastic, etc.) clearly should be considered. However, a general phenomenon emerges: increasing the number $m$ of lags while $n$ is fixed  helps (see Propositions~\ref{Prop.GaussAsymptotic} and~\ref{Prop.ranktest}) improve the chi-square approximation, but comes at the cost of additional degrees of freedom, hence larger critical values and lower rejection frequencies.  The choice of a test statistic, thus, relies on a trade-off:  once it is large enough for the level constraint to be satisfied,  increasing $m$  is not desirable unless a gain of power compensates for the loss caused by the additional degrees of freedom. 

\subsection{Robustness}

Next, we investigate the robustness properties of the center-outward rank-based and pseudo-Gaussian tests. A thorough theoretical treatment of this issue is beyond the scope of this paper, and we only perform a limited Monte-Carlo study,    focusing on two cases: resistance to a temporary change in the innovation density, with a short ``crisis period'' of heavy-tailed innovations,  and resistance to the presence of a patchy outlier. In both cases, the level and power of the  rank-based test remains largely unaffected while the pseudo-Gaussian  severely over-rejects.

\subsubsection{Resistance to  changes in the innovation density}\label{Sec.RobSkewt2}
In this section, we consider, for $-29\leq t\leq 1000$ (hence, $n=1030$ observations), the same VARMA(1,1) and VARMA(1,2) data-generating processes as in Section~\ref{Sec:sizepower}, with spherical Gaussian innovations for $t=1,\ldots,1000$ but, for the initial period~$t=-29,\ldots,0$ (that is, $3\%$ of the observation period), innovation outliers from a centered skew-$t_2$. 
This can be interpreted as a crisis-exit scenario with crisis period ending at time $t=0$---an information that is not available to the analyst. 


\begin{figure}[t!] 
\centering 
\includegraphics[width=130mm, height=64mm]{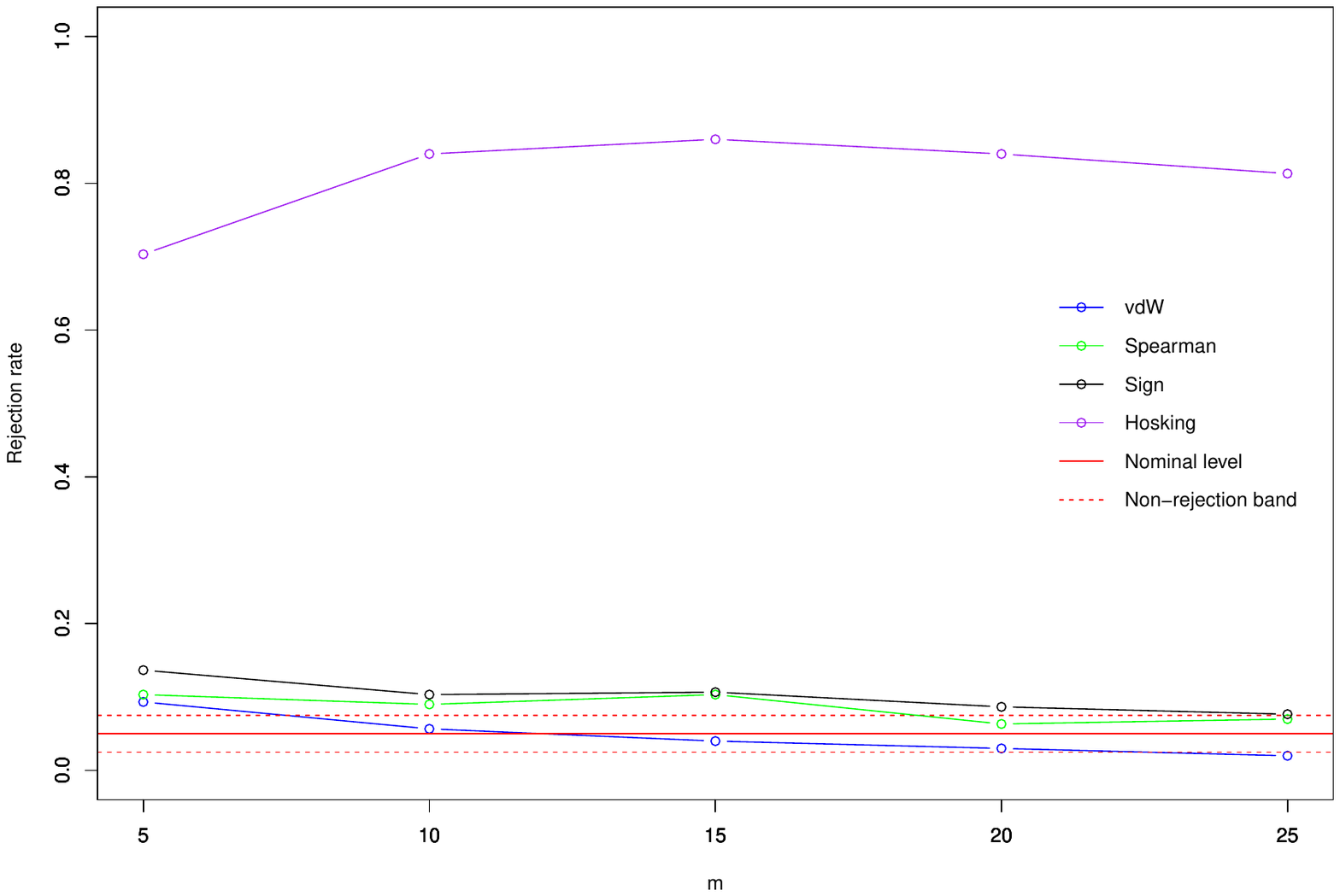}\vspace{-7mm}
\caption{\small Rejection frequencies (nominal level 5\%; asymptotic chi-square critical values), for~$m =~\!5,\, 10, \ldots, \,25$,  of the  Gaussian, van der Waerden, Spearman, and sign portmanteau tests for unspecified VARMA(1,1) model,  under the VARMA(1,1) model \eqref{VARMA(1,1)} with  spherical normal innovations  contaminated by $30$ points of skew-$t_2$ ones.   Number of observa\-tions~$n=1030$;  number of repli\-cations~$N=300$. The solid and dashed horizontal lines indicate the nominal level~$\alpha=5\%$ and the rejection limits of the~$5\%$ two-sided test  of the hypothesis that the actual level indeed is~$5\%$. }\label{plot.SizeAOskewt2}\vspace{3mm}
\end{figure}
 \begin{figure}[h!]
\centering 
\includegraphics[width=130mm, height=64mm]{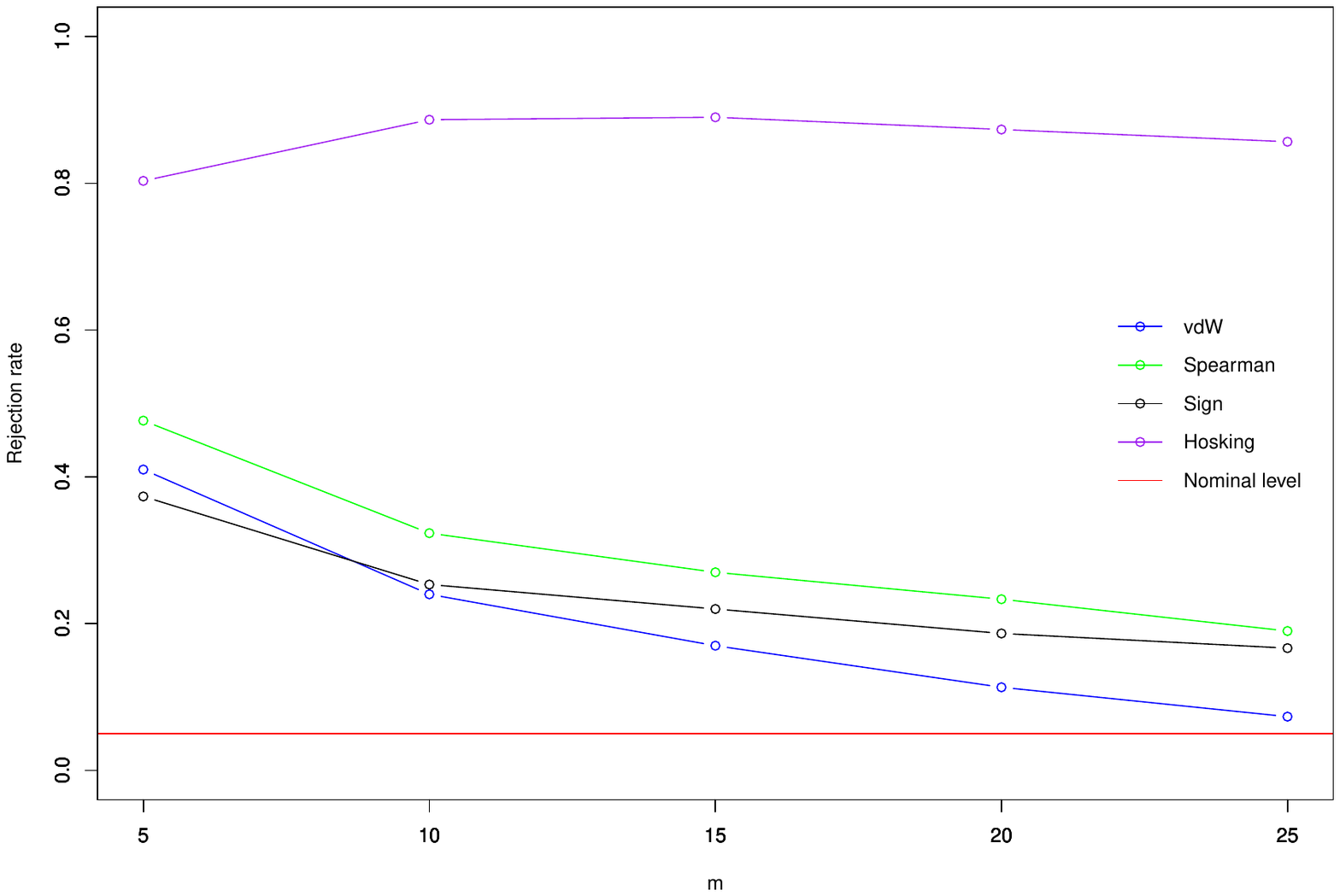}\vspace{-7mm}
\caption{\small Rejection frequencies (nominal level 5\%; asymptotic chi-square critical values),   for~$m\!=~\!5,\,10, \ldots,\, 25$, 
 of the  Gaussian, van der Waerden, Spearman, and sign portmanteau tests for unspecified VARMA(1,1) model, under the  VARMA(1,2) alternative \eqref{alternative}, with spherical normal innovations  contaminated by~$30$ points of skew-$t_2$  additive outliers.   Number of observations $n=1030$;  number of replica\-tions~\!$N=300$. The solid and dashed horizontal lines indicate the nominal level~$\alpha=5\%$. }\label{plot.PowerAOSkewt2}
\end{figure}

The rejection frequencies under the VARMA(1,1) null hypothesis \eqref{VARMA(1,1)} and  the VARMA(1,2) alternative \eqref{alternative} are shown in Figures~\ref{plot.SizeAOskewt2} and~\ref{plot.PowerAOSkewt2}, respectively. Comparing Figure~\ref{plot.SizeAOskewt2} and the top panel of Figure~\ref{plot.Size} reveals that the size of the rank-based tests under the null are remarkably robust (the van der Waerden satisfies the nominal $5\%$ constraint for $m\geq 5$, as it did in the upper panel of Figure~\ref{plot.Size}) while the  size of Hosking's  pseudo-Gaussian test, with a  rejection frequency under the null of about $0.80$, is very severely affected: the impact of the initial crisis ($3\%$ of the observation period), thus, is quite persistent.  The high rejection frequency of the classical portmanteau test, in Figure~\ref{plot.PowerAOSkewt2}, which is more or less the same as under the null, is meaningless 
and  entirely explained by overall over-rejection under the null. Quite on the contrary, the empirical powers of the rank-based tests remain quite stable and van der Waerden with $m=6$ remains the best choice (note that Spearman with~$m=14$ slightly over-rejects). Rank-based portmanteau tests, thus, unlike the traditional ones, are robust to the fact that 30 observations still belonging to the crisis period have been included in the analysis (the end-of-crisis date being unknown).

 \subsubsection{Resistance to  patchy outliers}
 
Another form of robustness is resistance to the presence of {\it patchy innovation outliers}. Here, at $t=500, 501$ and $t=502, 503$ respectively, shocks of sizes~$(20, 20)\pr$ and $(-20, -20)\pr$ are added to the Gaussian innovations for the same VARMA(1,1) and VARMA(1,2) data-generating processes as in Section~\ref{Sec:sizepower} (same observation period, $n=1000$). 

 \begin{figure}[b!] 
\centering 
\includegraphics[width=130mm, height=64mm]{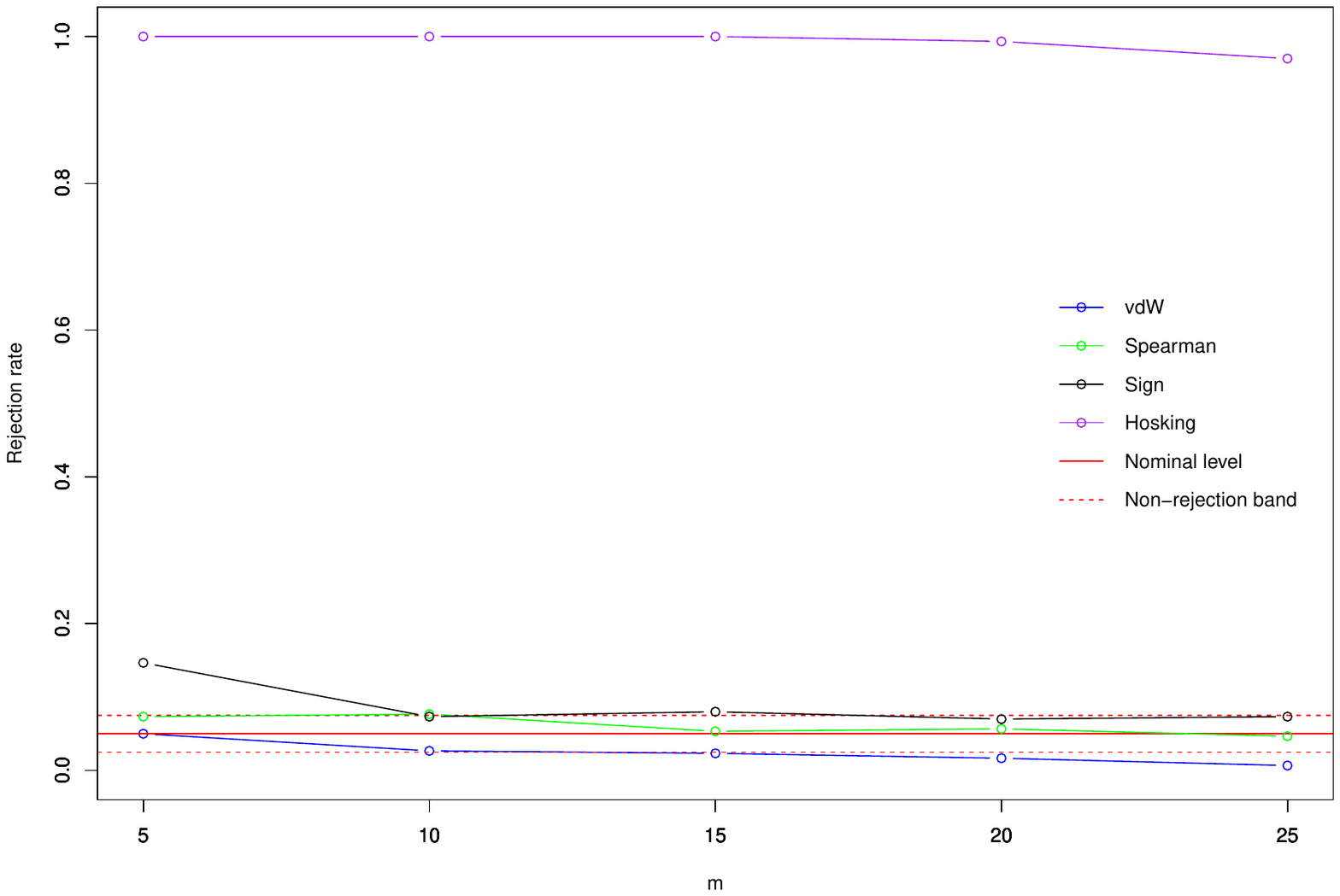}\vspace{-7mm}
\caption{\small Rejection frequencies (nominal level 5\%; asymptotic chi-square critical values), for~$m =~\!5,\, 10, \ldots, \,25$,  of the  Gaussian, van der Waerden, Spearman, and sign portmanteau tests for unspecified VARMA(1,1) model,  under the VARMA(1,1) model \eqref{VARMA(1,1)} with  spherical normal innovations  contaminated by $4$ points of additive outliers of size $20$.  Number of observations $n=1000$;  number of repli\-cations~$N=300$. The solid and dashed horizontal lines indicate the nominal level~$\alpha=5\%$ and the rejection limits of the~$5\%$ two-sided test  of the hypothesis that the actual level indeed is~$5\%$. }\label{plot.SizeAO4points}\vspace{-5mm}
\end{figure}
 \begin{figure}[t!]
\centering 
\includegraphics[width=130mm, height=64mm]{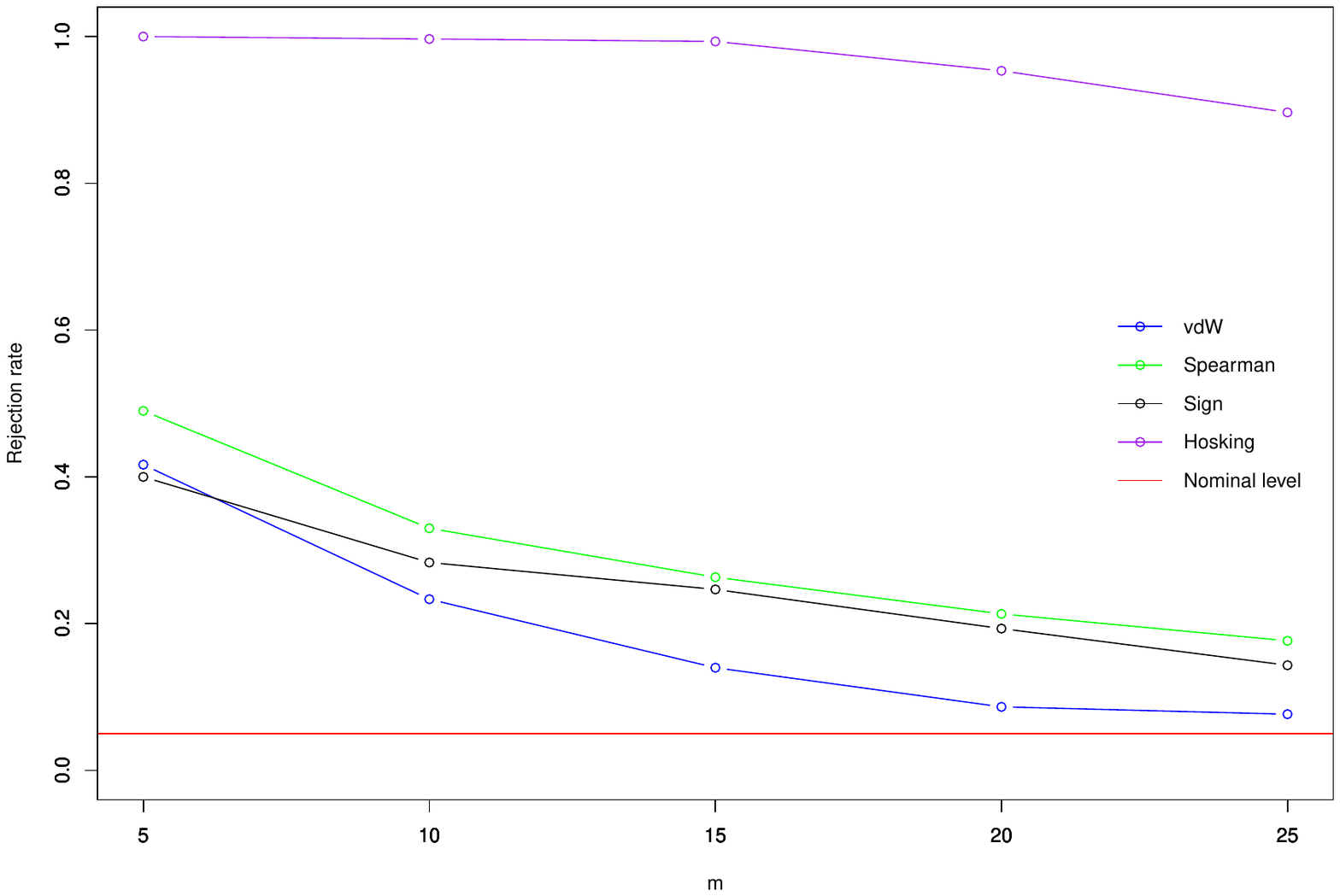}\vspace{-7mm}
\caption{\small Rejection frequencies (nominal level 5\%; asymptotic chi-square critical values),   for~$m\!=~\!5,\,10, \ldots,\, 25$, 
 of the  Gaussian, van der Waerden, Spearman, and sign portmanteau tests for unspecified VARMA(1,1) model, under the  VARMA(1,2) alternative \eqref{alternative}, with spherical normal innovations  contaminated by $4$ points of additive outliers of size $20$.   Number of observations $n=1030$;  number of replica\-tions~\!$N=300$. The solid and dashed horizontal lines indicate the nominal level~$\alpha=5\%$. }\label{plot.PowerAO4points}
\end{figure}

The rejection frequencies under the VARMA(1,1) null hypothesis and   VARMA(1,2) alternative  are shown in Figures~\ref{plot.SizeAO4points} and~\ref{plot.PowerAO4points}, respectively. 
 The plots clearly indicate that the rank-based tests are considerably more robust than their Gaussian competitor:  the empirical size of the Hosking  test (value~$1$ for $m = 5, 10,  15$, values  $0.993$ and $0.970$ for $m=20$ and $m=25$) is exploding and  its empirical power, which is close to the empirical size, is therefore meaningless. In stark contrast, Figure~\ref{plot.SizeAO4points} and the top panel of Figure~\ref{plot.Size} remain relatively similar:   the empirical sizes of the rank-based tests are not affected by outliers (e.g., for $m = 5$ or 6, the van der Waerden tests yield the nominal $5\%$ size) and their empirical powers remain  approximately the same as under the uncontaminated case.

\section{Conclusions} This paper achieves two objectives: a rigorous statement of the asymptotic behavior of the portmanteau test statistic under the null hypothesis, and the construction of rank-based versions of the same. Simulations indicate that the latter, thanks to a faster convergence, both in terms of the number $n$ of observations and the number $m$ of lags involved, bring substantial potential gains    of power when compared to their  classical counterparts---particularly so under skewed and heavy-tailed innovation densities. 

The same simulations also provides empirical evidence of a  better resistance of rank-based tests to the presence of contaminated innovations. A more thorough analysis, on the model of  \cite{HRRS86},  \cite{RY86}, or \cite{HR94} (all restricted to models with independent observations), of the robustness properties of rank tests in the time-series context is highly desirable, and should be left for further research.

\section*{Acknowledgement}
 The authors  thank the editor, associated editor, and three anonymous referees for helpful comments, which significantly helped improve the final version of the manuscript.

\singlespacing

\bigskip\bigskip

\newpage

\appendix

\noindent {\Large{\bf Appendix.}}\vspace{8mm}


\noindent \textbf{Proof of Proposition~\ref{asy.Gami.N}}. \bigskip

The result follows from deriving the asymptotic joint distribution, under~${\rm P}^{(n)}_{\tbth ;f}$,  of
$$(n-i)^{1/2} {\rm vec}({\bGamma}_{i; \mathcal{N}}^{(n)}(\bth)) ,\quad
(n-j)^{1/2} {\rm vec} ({\bGamma}_{j; \mathcal{N}}^{(n)}(\bth)),\quad\text{and}\quad \bDelta^{(n)}_{f} (\bth)$$ 
along the same lines as in the proof of Lemma~B.1 in \cite{HLL2022}. An application of Le Cam's third Lemma yields the asymptotic shifts in \eqref{aslinN} and concludes. 
 Details are left to the reader. \qed \vspace{5mm}

\noindent \textbf{Proof of Lemma~\ref{Lemport}}. \bigskip

Since $\bDelta\n_{\cal N}(\hat{\bth}\n_{\cal N}) = \0$,  we have, in view of~\eqref{=*},~${\bGamma}_{i; {\cal N}}^{(n)}(\hat{\bth}\n_{\cal N}) = {\bGamma}_{i; {\cal N}}^{(n)*}(\hat{\bth}\n_{\cal N})$. Moreover, \eqref{eq.linN} entails 
\begin{align*}
&(n-i)^{1/2} {\rm vec} ({\bGamma}_{i; {\cal N}}^{(n)*}(\hat{\bth}\n_{\cal N}) - {\bGamma}_{i; {\cal N}}^{(n)*}(\bth))   \\
&= (n\! -\! i)^{1/2} {\rm vec} ({\bGamma}^{(n)}_{i; {\cal N}}(\hat{\bth}\n_{\cal N})\! -\! {\bGamma}^{(n)}_{i; {\cal N}}(\bth)) 
\! -\! ({\bf \Sigma}\otimes {\bf \Sigma}^{-1}) \c_{i,\tbth}\pr\Big( \sum_{i=1}^{n-1}\! \c_{i,\tbth} ({\bf \Sigma}\otimes {\bf \Sigma}^{-1}) \c_{i,\tbth}\pr\!\Big)^{-1}  \\
&\ \hspace{35mm}  \times\left(\sum_{i=1}^{n-1} \c_{i, \tbth}  (n-i)^{1/2} {\rm vec}({\bGamma}^{(n)}_{i; {\cal N}}(\hat{\bth}\n_{\cal N}) - {\bGamma}^{(n)}_{i; {\cal N}}(\bth))\right)  + o_{\rm P}(1) \\
&= \left( ({\bf \Sigma}\otimes {\bf \Sigma}^{-1}) \c_{i, \tbth}\pr - ({\bf \Sigma}\otimes {\bf \Sigma}^{-1}) \c_{i, \tbth}\pr  \Big( \sum_{i=1}^{n-1} \c_{i,\tbth} ({\bf \Sigma}\otimes {\bf \Sigma}^{-1}) \c_{i,\tbth}\pr\Big)^{-1} \right.  \\
&\left.\ \qquad\qquad\qquad\qquad\qquad\quad \times \Big( \sum_{i=1}^{n-1} \c_{i,\tbth} ({\bf \Sigma}\otimes {\bf \Sigma}^{-1}) \c_{i,\tbth}\pr\Big) \right) n^{1/2}(\hat{\bth}\n_{\cal N} - \bth)  + o_{\rm P}(1) \\
&= o_{\rm P}(1).
\end{align*}
The result follows. \qed \vspace{5mm}

\noindent \textbf{Proof of Proposition~\ref{Prop.GaussAsymptotic}}. \bigskip

Due to Lemma~\ref{Lemport} and the exponential decrease, as $i\to\infty$,  of $\Vert\c_{i, \tbth}\Vert$, $$\Vert (n-i)^{1/2} {\rm vec} \left({\bGamma}^{(m, n)**}_{i, {\cal N}}(\bth) - {\bGamma}_{i, {\cal N}}^{(n)}(\hat{\bth}^{(n)}_{{\cal N}}) \right) \Vert$$
with probability one  converges to zero exponentially fast as $m$ increases and~$n\rightarrow \infty$. Part {\it (i)} then follows.

Turning to Part {\it (ii)}, let
\begin{align*}
{\bGamma}_{{\cal N}}^{(m, n)**}(\bth) \coloneqq &   n^{-1/2} \left( (n-1)^{1/2} (\text{vec}(  {\bGamma}_{1; {\cal N}}^{(m, n)**}(\bth)))^\prime, \ldots\right. \\    
&\hspace{40mm}\left.  \ldots, (n-m)^{1/2} (\text{vec}(  {\bGamma}_{m; {\cal N}}^{(m, n)**}(\bth)))^\prime \right)^\prime .
\end{align*}
Note that ${\bGamma}_{{\cal N}}^{(m, n)**}(\bth)$ can be written as
${\bGamma}_{{\cal N}}^{(m, n)**}(\bth) = \E^{(m)}_{{\cal N}}(\bth) {\bGamma}_{{\cal N}}^{(m, n)}(\bth),$
where 
\begin{align*}
\E^{(m)}_{{\cal N}}(\bth) &\coloneqq  \I_{md^2}
 - (\I_{m} \otimes {\bf \Sigma}\otimes {\bf \Sigma}^{-1}) \C_{\tbth}^{(m+1)\prime}   \\  
 &\hspace{22mm}
\times  \left( \C_{\tbth}^{(m+1)} (\I_{m} \otimes {\bf \Sigma}\otimes {\bf \Sigma}^{-1}) \C_{\tbth}^{(m+1)\prime}
 \right)^{-1} \C_{\tbth}^{(m+1)}
 \end{align*}
is an idempotent matrix. Then it follows from Proposition~\ref{asy.Gami.N} that ${\bGamma}_{{\cal N}}^{(m, n)**}(\bth)$ is asymptotically normal with mean zero and covariance $\E^{(m)}_{{\cal N}}(\bth) (\I_{m} \otimes {\bf \Sigma} \otimes {\bf \Sigma})$ under ${\rm P}^{(n)}_{\tbth; f}$. It remains to prove that 
 ${\rm tr}(\E^{(m)}_{\tbth}) = (m-p-q)d^2,$ 
where ${\rm tr}(\E^{(m)}_{\tbth})$ denotes the trace of $\E^{(m)}_{\tbth}$. Using ${\rm tr}(\A\B\C) = {\rm tr}(\C\A\B)$, we obtain
\begin{align*}
{\rm tr}(\E^{(m)}_{\tbth}) &= {\rm tr}(\I_{md^2}) - {\rm tr}\left( \C_{\tbth}^{(m+1)} (\I_{m} \otimes {\bf \Sigma}\otimes {\bf \Sigma}^{-1}) \C_{\tbth}^{(m+1)\prime} \right.\\
&\hspace{50mm} \left. \times \left( \C_{\tbth}^{(m+1)} (\I_{m} \otimes {\bf \Sigma}\otimes {\bf \Sigma}^{-1}) \C_{\tbth}^{(m+1)\prime}\right)^{-1} \right) \\
&= md^2 - {\rm tr}(\I_{(p+q)d^2})  = (m-p-q)d^2.
\end{align*}
The result follows. \qed \vspace{5mm}


\noindent \textbf{Proof of Lemma~\ref{Lemport.rank}}. \bigskip

It follows from \eqref{asylin2} and Lemma~4.4 in \cite{Kreiss87} that, under  ${\rm P}^{(n)}_{\tbth; f}$, 
\[\tenq{\bDelta}^{(n)}_{\J_1, \J_2}(\tenq{\bth}^{(n)}_{\J_1, \J_2}) = \tenq{\bDelta}^{(n)}_{\J_1, \J_2}(\bth) - n^{1/2} {\bUpsilon}\n_{\J_1, \J_2, f}(\bth)   \left(
\tenq{\bth}^{(n)}_{\J_1, \J_2}-\bth \right) + o_{\rm P}(1).
\]
Hence, letting~${\bUpsilon}_{\J_1, \J_2, f}(\bth)\coloneqq  \lim_{n\to\infty} {\bUpsilon}\n_{\J_1, \J_2, f}(\bth)$, it follows from the definition of $\tenq{\bth}^{(n)}_{\J_1, \J_2}$ in \eqref{def.Rest} that 
\begin{align*}
\tenq{\bDelta}^{(n)}_{\J_1, \J_2}(\tenq{\bth}^{(n)}_{\J_1, \J_2})=\,& \tenq{\bDelta}^{(n)}_{\J_1, \J_2}(\bth) -
n^{1/2}\left( \bUpsilon_{\J_1, \J_2, f}(\bth) \bar{\bth}\n \right. \\
&\left. \hspace{30mm}+ n^{-1/2}\tenq{\bDelta}^{(n)}_{\J_1, \J_2}( \bar{\bth}\n) -  {\bUpsilon}_{\J_1, \J_2, f}(\bth)\bth
\right) 
+ o_{\rm P}(1)\\ 
=\,& \tenq{\bDelta}^{(n)}_{\J_1, \J_2}(\bth) -n^{1/2}{\bUpsilon}_{\J_1, \J_2, f}(\bth)\left(\bar{\bth}\n -\bth
\right) -\tenq{\bDelta}^{(n)}_{\J_1, \J_2}( \bar{\bth}\n) + o_{\rm P}(1)\\ 
=\,&n^{1/2} {\bUpsilon}_{\J_1, \J_2, f}(\bth)\left(\bar{\bth}\n -\bth
\right) -{n^{1/2}\bUpsilon}_{I; \J_1, \J_2, f}(\bth)\left(\bar{\bth}\n -\bth
\right)+ o_{\rm P}(1)\\ =\,& o_{\rm P}(1).
\end{align*}
This establishes part {\it (i)} of the lemma. Part {\it (ii)} follows as a corollary of Part~{\it (i)}, since $\tenq{\bGamma}_{i, \J_1, \J_2, f}^{(n)*}(\tenq{\bth}^{(n)}_{\J_1, \J_2})$ is~the residual of the regression of $\tenq{\bGamma}_{i, \J_1, \J_2}^{(n)}(\bth)$ with respect to $\tenq{\bDelta}^{(n)}_{\J_1, \J_2}(\bth)$ computed at~$\bth=\tenq{\bth}^{(n)}_{\J_1, \J_2}$.  Asymptotic linearity, the asymptotic discreteness of  $\tenq{\bth}^{(n)}_{\J_1, \J_2}$ and Lemma~4.4 of~\cite{Kreiss87} entail~{\it (iii)}.\qed \vspace{5mm}

\noindent \textbf{Proof of Proposition~\ref{Prop.ranktest}}. \bigskip

Part {\it (i)} follows from~\eqref{eq.PPstar} and the exponential decrease of $\Vert\c_{i, \tbth}\Vert$. Turning to Part {\it (ii)}, let
\begin{align*}
\tenq{\bGamma}_{\J_1, \J_2, f}^{(m, n)**}(\bth) \coloneqq &   n^{-1/2} \left( (n-1)^{1/2} (\text{vec}(  \tenq{\bGamma}_{1, \J_1, \J_2, f}^{(m, n)**}(\bth)))^\prime, \ldots \right. \\
&\hspace{40mm}\left. \ldots,   (n-m)^{1/2} (\text{vec}(  \tenq{\bGamma}_{m, \J_1, \J_2, f}^{(m, n)**}(\bth)))^\prime \right)^\prime
\end{align*}
where $\tenq{\bGamma}^{(m, n)**}_{i, \J_1, \J_2, f}(\bth)$ is defined in~\eqref{Gamstarstar}. Note that $\tenq{\bGamma}_{\J_1, \J_2, f}^{(m, n)**}(\bth)$ can be written as
$$\tenq{\bGamma}_{\J_1, \J_2, f}^{(m, n)**}(\bth) = \E^{(m)}_{\J_1, \J_2, f}(\bth) \tenq{\bGamma}_{\J_1, \J_2}^{(m, n)}(\bth).$$
Then,  from the definition of $\tenq{Q}_{m; \J_1, \J_2, f}^{(n)**} (\bth)$ in~\eqref{Pstarstar}, 
\begin{align*}
\tenq{Q}_{m; \J_1, \J_2, f}^{(n)**} (\bth) &= n (\tenq{\bGamma}_{\J_1, \J_2, f}^{(m, n)**}(\bth))\pr \left({\rm diag}({\bOmega}^{(m)}_{i, \J_1, \J_2; f}(\bth))_{1\leq i \leq m}\right)^{-} \tenq{\bGamma}_{\J_1, \J_2, f}^{(m, n)**}(\bth) \\
&= n (\tenq{\bGamma}_{\J_1, \J_2}^{(m, n)}(\bth))\pr \E^{(m)\prime}_{\J_1, \J_2, f}(\bth) \left({\rm diag}({\bOmega}^{(m)}_{i, \J_1, \J_2; f}(\bth))_{1\leq i \leq m}\right)^{-}  \\ 
&\hspace{60mm} \times \E^{(m)}_{\J_1, \J_2, f}(\bth) \tenq{\bGamma}_{\J_1, \J_2}^{(m, n)}(\bth)\\&= n (\tenq{\bGamma}_{\J_1, \J_2}^{(m, n)}(\bth))\pr \E^{(m)\prime}_{\J_1, \J_2, f}(\bth) \left(\E^{(m)}_{\J_1, \J_2, f}(\bth) (\I_m \otimes \D_{\J_1, \J_2})\right. \\
& \hspace{50mm}\left. \times  \E^{(m)\prime}_{\J_1, \J_2, f}(\bth) \right)^{-}  \E^{(m)}_{\J_1, \J_2, f}(\bth) \tenq{\bGamma}_{\J_1, \J_2}^{(m, n)}(\bth),
\end{align*}
where the last equality follows from \eqref{eq.diagOmega}. 
Now, letting
\begin{align*}
\M_{\tbth} &\coloneqq \M^{(m)}_{\J_1, \J_2, f}(\bth)
\\ 
& \coloneqq (\I_m \otimes \D_{\J_1, \J_2}^{1/2}) \E^{(m)\prime}_{\J_1, \J_2, f}(\bth) \left(\E^{(m)}_{\J_1, \J_2, f}(\bth) (\I_m \otimes \D_{\J_1, \J_2}) \E^{(m)\prime}_{\J_1, \J_2, f}(\bth) \right)^{-} \\
&\hspace{70mm} \times \E^{(m)}_{\J_1, \J_2, f}(\bth) (\I_m \otimes \D_{\J_1, \J_2}^{1/2}),
\end{align*}
we obtain 
$$\tenq{Q}_{m; \J_1, \J_2, f}^{(n)**} (\bth) = n (\tenq{\bGamma}_{\J_1, \J_2}^{(m, n)}(\bth))\pr (\I_m \otimes \D_{\J_1, \J_2}^{-1/2}) \M_{\tbth} (\I_m \otimes \D_{\J_1, \J_2}^{-1/2})  \tenq{\bGamma}_{\J_1, \J_2}^{(m, n)}(\bth).$$
In view of Proposition~\ref{asy.Gami.utilde}, 
$(\I_m \otimes \D_{\J_1, \J_2}^{-1/2})  \tenq{\bGamma}_{\J_1, \J_2}^{(m, n)}(\bth)$ is asymptotically  normal with covariance $\I_{md^2}$. Moreover, since $\M_{\tbth}$ is symmetric and idempotent~(indeed, $\M_{\tbth} \M_{\tbth} = \M_{\tbth}$), it remains to show that $\M_{\tbth}$ is of rank $(m-p-q)d^2$. Note that $\M_{\tbth}$ has the same rank as $\E^{(m)}_{\J_1, \J_2, f}(\bth)$ which is also idempotent, 
it suffices to prove that 
 ${\rm tr}(\E^{(m)}_{\J_1, \J_2, f}(\bth)) = (m-p-q)d^2.$ 
Using again the fact that ${\rm tr}(\A\B\C) = {\rm tr}(\C\A\B)$, we have
\begin{align*}
{\rm tr}&(\E^{(m)}_{\J_1, \J_2, f}(\bth))  \\ 
&= {\rm tr}(\I_{md^2}) - {\rm tr}((\I_m \otimes  \K_{\J_1, \J_2, f}) \C^{(m+1)\prime}_{\tbth} (\sum_{i=1}^{m} \c_{i, \tbth} \K_{\J_1, \J_2, f} \c^{\prime}_{i, \tbth})^{-1} \C^{(m+1)}_{\tbth}) \\
&= md^2 - {\rm tr}(\C^{(m+1)}_{\tbth}(\I_m \otimes  \K_{\J_1, \J_2, f}) \C^{(m+1)\prime}_{\tbth} (\sum_{i=1}^{m} \c_{i, \tbth} \K_{\J_1, \J_2, f} \c^{\prime}_{i, \tbth})^{-1}) \\
&= md^2 - {\rm tr}(\I_{(p+q)d^2})  = (m-p-q)d^2.
\end{align*}
The result follows. \qed \vspace{3mm}

\end{document}